%% file: rf.tex
\documentclass[twoside
]{amsart}%

  \usepackage{a4wide}
  \usepackage{amsmath}
  \usepackage{amssymb}
  \usepackage{amsthm}
  \usepackage{mathrsfs}
  \usepackage{bbm}
  \usepackage{pifont}
  \usepackage{color}
  \usepackage{graphicx}

  \usepackage[normalem]{ulem}
	
	\def\AA{{\ifmmode{\mathbbm{A}}\else{$\mathbbm{A}$}\fi}}
	\def\BB{{\ifmmode{\mathbbm{B}}\else{$\mathbbm{B}$}\fi}}
	\def\CC{{\ifmmode{\mathbbm{C}}\else{$\mathbbm{C}$}\fi}}
	\def\EE{{\ifmmode{\mathbbm{E}}\else{$\mathbbm{E}$}\fi}}
	\def\FF{{\ifmmode{\mathbbm{F}}\else{$\mathbbm{F}$}\fi}}
	\def\HH{{\ifmmode{\mathbbm{H}}\else{$\mathbbm{H}$}\fi}}
	\def\KK{{\ifmmode{\mathbbm{K}}\else{$\mathbbm{K}$}\fi}}
	\def\NN{{\ifmmode{\mathbbm{N}}\else{$\mathbbm{N}$}\fi}}
	\def\PP{{\ifmmode{\mathbbm{P}}\else{$\mathbbm{P}$}\fi}}
	\def\QQ{{\ifmmode{\mathbbm{Q}}\else{$\mathbbm{Q}$}\fi}}
	\def\RR{{\ifmmode{\mathbbm{R}}\else{$\mathbbm{R}$}\fi}}
	\def\TT{{\ifmmode{\mathbbm{T}}\else{$\mathbbm{T}$}\fi}}
	\def\UU{{\ifmmode{\mathbbm{U}}\else{$\mathbbm{U}$}\fi}}
	\def\ZZ{{\ifmmode{\mathbbm{Z}}\else{$\mathbbm{Z}$}\fi}}
	
	\def\A{{\ifmmode{\mathscr{A}}\else{$\mathscr{A}$}\fi}}
	\def\B{{\ifmmode{\mathscr{B}}\else{$\mathscr{B}$}\fi}}
	\def\C{{\ifmmode{\mathscr{C}}\else{$\mathscr{C}$}\fi}}
	\def\D{{\ifmmode{\mathscr{D}}\else{$\mathscr{D}$}\fi}}
	\def\E{{\ifmmode{\mathscr{E}}\else{$\mathscr{E}$}\fi}}
	\def\F{{\ifmmode{\mathscr{F}}\else{$\mathscr{F}$}\fi}}
	\def\G{{\ifmmode{\mathscr{G}}\else{$\mathscr{G}$}\fi}}
	\def\H{{\ifmmode{\mathscr{H}}\else{$\mathscr{H}$}\fi}}
	\def\I{{\ifmmode{\mathscr{I}}\else{$\mathscr{I}$}\fi}}
	\def\J{{\ifmmode{\mathscr{J}}\else{$\mathscr{J}$}\fi}}
	\def\K{{\ifmmode{\mathscr{K}}\else{$\mathscr{K}$}\fi}}
	\def\L{{\ifmmode{\mathscr{L}}\else{$\mathscr{L}$}\fi}}
	\def\M{{\ifmmode{\mathscr{M}}\else{$\mathscr{M}$}\fi}}
	\def\N{{\ifmmode{\mathscr{N}}\else{$\mathscr{N}$}\fi}}
	\def\O{{\ifmmode{\mathscr{O}}\else{$\mathscr{O}$}\fi}}
	\def\P{{\ifmmode{\mathscr{P}}\else{$\mathscr{P}$}\fi}}
	\def\Q{{\ifmmode{\mathscr{Q}}\else{$\mathscr{Q}$}\fi}}
	\def\R{{\ifmmode{\mathscr{R}}\else{$\mathscr{R}$}\fi}}
	\def\S{{\ifmmode{\mathscr{S}}\else{$\mathscr{S}$}\fi}}
	\def\T{{\ifmmode{\mathscr{T}}\else{$\mathscr{T}$}\fi}}
	\def\U{{\ifmmode{\mathscr{U}}\else{$\mathscr{U}$}\fi}}
	\def\V{{\ifmmode{\mathscr{V}}\else{$\mathscr{V}$}\fi}}
	\def\W{{\ifmmode{\mathscr{W}}\else{$\mathscr{W}$}\fi}}
	\def\X{{\ifmmode{\mathscr{X}}\else{$\mathscr{X}$}\fi}}
	\def\Y{{\ifmmode{\mathscr{Y}}\else{$\mathscr{Y}$}\fi}}
	\def\Z{{\ifmmode{\mathscr{Z}}\else{$\mathscr{Z}$}\fi}}

	\newtheoremstyle{slanted}
	{}
	{}
	{\slshape}
	{}
	{\bfseries}
	{.}
	{ }
	{}
	
	\theoremstyle{slanted}
	\newtheorem{theo}{Theorem}[section]
	\newtheorem{prop}[theo]{Proposition}
	
	\newtheorem{remark}[theo]{Remark}
	
	\newtheorem{lemma}[theo]{Lemma}
	\newtheorem{definition}[theo]{Definition}
	\newtheorem{corollary}[theo]{Corollary}
	
	\def\ind#1{\mathbbmss{1}_{#1}}
	\def\egdef{:=}

	\newcommand{\tend}[2]{\xrightarrow[#1\to#2]{}}

	\def\ind#1{\mathbbmss{1}_{#1}}

\title{How do random Fibonacci sequences grow?}

\author{\'Elise Janvresse, Beno\^it Rittaud, Thierry de la Rue}

\address{\'Elise Janvresse, Thierry de la Rue:
Laboratoire de Math\'ematiques Rapha\"el Salem, 
Universit\'e de Rouen, CNRS -- 
Avenue de l'Universit\'e -- 
F76801 Saint \'Etienne du Rouvray.}
\email{Elise.Janvresse@univ-rouen.fr\\Thierry.de-la-Rue@univ-rouen.fr}
\address{Beno\^it Rittaud: Laboratoire Analyse, G\'eom\'etrie et Applications, Universit\'e Paris 13 Institut Galil\'ee, CNRS -- 
99 avenue Jean-Baptiste Cl\'ement -- 
F93 430 Villetaneuse.}
\email{rittaud@math.univ-paris13.fr}
\begin{document}
\bibliographystyle{amsplain}

\begin{abstract}
We study the random Fibonacci sequences defined by 
$F_1=F_2=\widetilde F_1=\widetilde F_2=1$ and for $n\ge 1$, 
$F_{n+2} = F_{n+1} \pm F_{n}$ (linear case) and 
$\widetilde F_{n+2} = |\widetilde F_{n+1} \pm \widetilde F_{n}|$ (non-linear case), 
where each $\pm$ sign is independent and either $+$ with probability $p$ or $-$ with probability $1-p$ ($0<p\le 1$). 
Our main result is that the exponential growth of $F_n$ for $0<p\le 1$, and of $\widetilde F_{n}$
for $1/3\le p\le 1$ is almost surely given by
$$
\int_0^\infty \log x\, d\nu_\alpha (x), 
$$
where $\alpha$ is an explicit function of $p$ depending on the case we consider, and $\nu_\alpha$ is an explicit probability distribution on $\RR_+$ defined inductively on Stern-Brocot intervals. 

In the non-linear case, the largest Lyapunov exponent is not an analytic function of $p$, since we prove that it is equal to zero for $0<p\le1/3$.
We also give some results about the variations of the largest Lyapunov exponent, and provide a formula for its derivative.
\end{abstract}

\maketitle
\section{Introduction}
In this article, we wish to investigate the exponential growth of \emph{random Fibonacci sequences} $(F_n)_{n\ge 1}$ and $(\widetilde F_n)_{n\ge 1}$, defined inductively 
by $F_1=F_2=\widetilde F_1=\widetilde F_2=1$, and for all $n\ge 1$, 
\begin{equation}
\label{linear case}
F_{n+2} = F_{n+1} \pm F_{n} \qquad \mbox{(linear case)},
\end{equation}
\begin{equation}
\label{non-linear case}
\widetilde F_{n+2} = |\widetilde F_{n+1} \pm \widetilde F_{n}| \qquad \mbox{(non-linear case)}, 
\end{equation}
where each $\pm$ sign is independent and either $+$ with probability $p$ or $-$ with probability $1-p$ ($0<p\le 1$). 
In the case $p=1/2$, $(|F_n|)$ and $(\widetilde F_n)$ have the same distribution law as the sequence $(|t_n|)$ studied by Viswanath \cite{viswanath2000}. In his paper, using Furstenberg's formula \cite{furstenberg1963} (see also \cite{bougerol1985}, Chapter~II), Viswanath proves that with probability 1,
$$ \sqrt[n]{|t_n|}\tend{n}{+\infty} 1.13198824\ldots,$$
where the logarithm of the number is computed as the integral of the function 
$$ m\longmapsto \dfrac{1}{4}\log\left(\dfrac{1+4m^4}{(1+m^2)^2}\right) $$
with respect to some explicit ``fractal'' measure $\nu_f$. 

Our purpose here is to give a formula for any parameter $p\in]0,1]$, and to provide some results on the dependence on $p$ of the upper Lyapunov exponents. 
By contrast with the linear case \eqref{linear case}, the non-linear case \eqref{non-linear case} cannot be viewed as a product of i.i.d. random matrices. This explains why the upper Lyapunov exponent in the non-linear case is not an analytic function of $p$, as can be seen in Theorem~\ref{MainTheorem}. 
Our method does not make use of Furstenberg's formula, but relies on the reduction of random Fibonacci sequences exposed in \cite{rittaud2006}.

\subsection{Results}
\label{res}
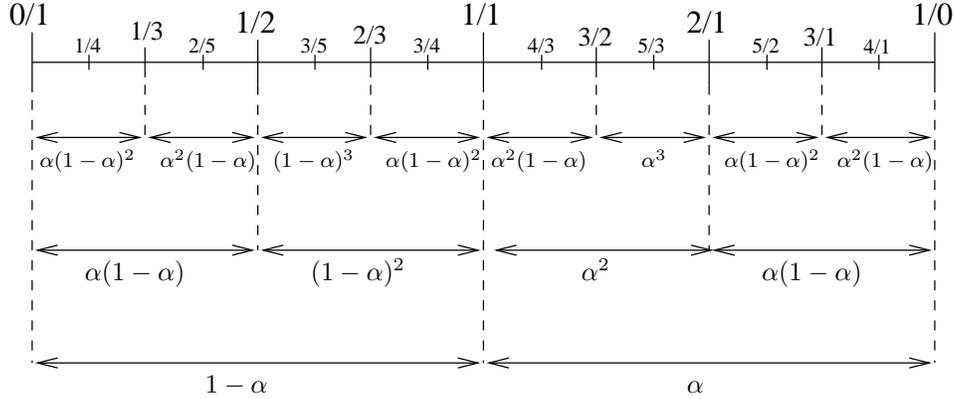
\begin{figure}[h]
	\begin{center}
	\input{mesure.pstex_t}
	\end{center}
\caption{The measure $\nu_\alpha$ on Stern-Brocot intervals of rank 1, 2, 3.
First assign mass $1-\alpha$ to $[0,1]$ and $\alpha$ to $[1,\infty]$. Once $\nu_\alpha$ is defined on some Stern-Brocot interval of rank $r$, a proportion $\alpha$ of its mass is given to the left (respectively right) subinterval of rank $r+1$ when $r$ is odd (respectively even). 
}
\label{mesure}
\end{figure}

Our main result is the following.
\begin{theo}
\label{MainTheorem}
\quad
\begin{itemize}
\item {\bf Linear case:}
$$ 
\dfrac{1}{n} \log |F_n| \tend{n}{\infty} \gamma_p = \int_0^\infty \log x\, d\nu_\alpha (x)\quad\mbox{a.s.}, $$
where 
$$
 \alpha \egdef \frac{3p-2+\sqrt{5p^2-8p+4}}{2p}\, .
$$
\item {\bf Non-linear case:} If $p\le 1/3$, $\dfrac{1}{n} \log \widetilde F_n \tend{n}{\infty} \widetilde \gamma_p = 0$. 
If $p\ge 1/3$, 
$$ 
\dfrac{1}{n} \log \widetilde F_n \tend{n}{\infty} \widetilde \gamma_p = \int_0^\infty \log x\, d\nu_\alpha (x)\quad\mbox{a.s.}, $$
where 
$$
 \alpha \egdef \frac{2p}{p + \sqrt{p(4-3p)}}\, .
$$
\end{itemize}
In both cases, $\nu_\alpha$ is an explicit probability distribution on $\RR_+$ defined inductively on Stern-Brocot intervals (see Section~\ref{section_distributionQ} and Figure~\ref{mesure}).
\end{theo}

For $p=1/2$, we get $\alpha=\phi^{-1}$ both in the linear and the non-linear cases, where $\phi\egdef (1+\sqrt{5})/2$ is the golden ratio. The measure $\nu_\alpha$ is then equal to Viswanath's fractal measure conditioned on $\RR^+$. 
For $p=1$, which corresponds to the classical Fibonacci sequence, $\alpha=1$ and $\nu_\alpha$ is the Dirac mass on $\phi$. 
When $p\to 0$ in the linear case, or $p\to 1/3$ in the non-linear case, $\alpha\to 1/2$ and $\nu_\alpha\to\nu_{1/2}$ which is the probability measure on $\RR^+$ giving the same mass $2^{-r}$ to each Stern-Brocot interval of rank $r$. This measure is related to Minkowski's Question Mark Function (see \cite{denjoy1938}): 
$$ 
\forall x\in[0,1],\quad ?(x) = 2\,\nu_{1/2}([0,x]).
$$

\begin{remark}
\label{moyenne}
The exponents $\gamma_p$ and ${\widetilde\gamma}_{p}$ correspond to 
almost-sure exponential growth of random Fibonacci sequences. 
We could also consider the average point of view, that is, ask for the limit of
${\dfrac{1}{n}\log(\EE(|{F}_n|))}$ and ${\dfrac{1}{n}\log(\EE({\widetilde F}_n))}$ 
(where $\EE$ stands  for the expectation). 

In the non-linear case, we know how to give an explicit expression of
the limit for any $p$. 
(Of course, by Jensen's inequality, this limit is bounded below by ${\widetilde \gamma}_p$.) 
It turns out that the critical value of $p$ for which this limit vanishes is $p=1/4$ (compare it with the value $1/3$ in Theorem~\ref{MainTheorem}). 
The techniques used to obtain this result are quite different from those presented
here; They mainly rely on ideas introduced in~\cite{rittaud2006}, in which the case
$p=1/2$ is treated. 
Details and proofs will be given in a forthcoming paper.
\end{remark}

\medskip
In Section~\ref{study}, we study some properties of the functions $p\mapsto\gamma_p$ and $p\mapsto \widetilde \gamma_p$, and prove the following theorem.

\begin{theo}
\label{increasing_gamma}
\quad
\begin{itemize}
\item {\bf Linear case:} $p\mapsto \gamma_p$ is an increasing function of $p$, satisfying
$$
\lim_{p\to 0} \gamma_p = 0, \quad \gamma_1 = \log\phi \quad \mbox{ and }\quad \frac{d\gamma_p}{dp}(1) = \frac{\log 5}{2}. 
$$
\item {\bf Non-linear case:} $p\mapsto \widetilde \gamma_p$ is a continuous function of $p$, increasing on $]1/3,1]$, satisfying
$$
\widetilde \gamma_{p} = 0 \mbox{ for } p\le 1/3, \quad 
\widetilde \gamma_1 = \log\phi \quad \mbox{ and }\quad 
\frac{d\widetilde \gamma_p}{dp}(1) = \frac{\log 5}{2}.
$$
\end{itemize}
\end{theo}
One of the ingredients for the proof is a formula for the derivative of $\gamma_p$ (or $\widetilde \gamma_{p}$) with respect to $\alpha$, involving the product measure $\nu_\alpha\otimes\nu_\alpha$ (see Proposition~\ref{prop_derivee_gamma}).

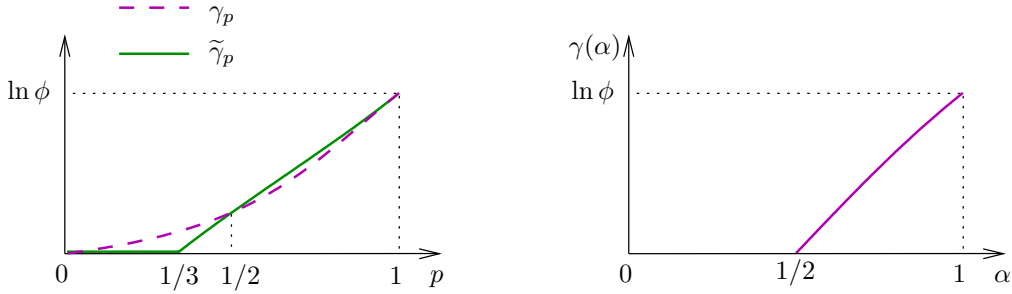
\begin{figure}[h]
	\begin{center}
	\input{gammas.pstex_t}
	\end{center}
\caption{Graphs of $\gamma_p$ and $\widetilde \gamma_{p}$ as a function of $p$ (left) and as a function of $\alpha$ (right). 
Notice that $\widetilde \gamma_{p} = 0$ for $p\in]0, 1/3]$. For $p\in[1/3, 1]$, the graph of $\widetilde \gamma_{p}$ is not the straight line it looks like: Compare the average slope with the derivative in $1$.}
\label{graphs_of_gamma}
\end{figure}

\section{Reduced sequences}
\label{section2}

\subsection{Random paths in $\mathbf{T}$ and $\mathbf{\widetilde T}$}
\label{subsection2.1}
\begin{figure}
	\begin{center}
	\includegraphics{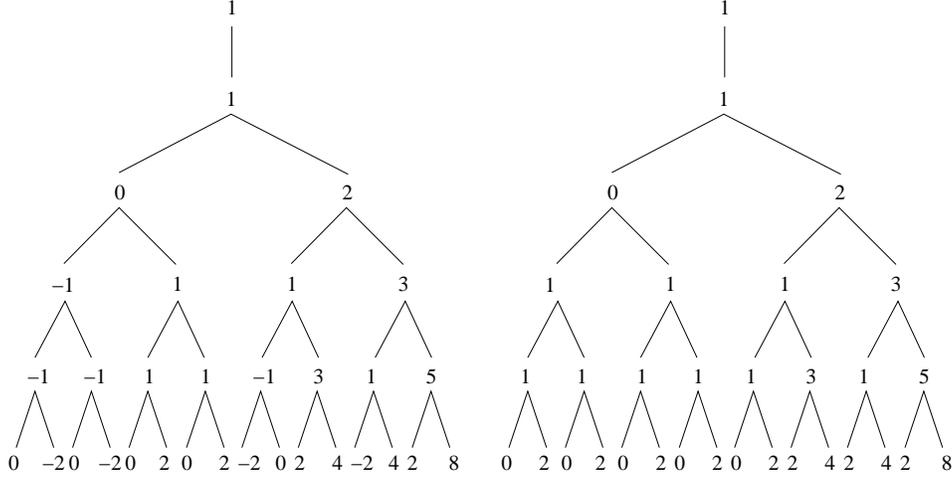}
	\end{center}
\caption{First lines of the trees $\mathbf{T}$ (left) and $\mathbf{\widetilde T}$ (right).}
\label{treeT}
\end{figure}

The sequences $(F_n)_n$ and $(\widetilde F_n)_n$ can be read along random paths in the trees $\mathbf{T}$ and $\mathbf{\widetilde T}$ described as follows. These two trees have the same structure, but differ by the labels attached to the vertices. 
Each vertex is labelled by an integer: The root and its only child are labelled by 1. Any other vertex has two children, left and right. If $v$ is a right child, its label is the sum of its father's and grandfather's labels; 
If $v$ is a left child, its label is the difference between its father's and grandfather's labels in the tree $\mathbf{T}$, whereas in the tree $\mathbf{\widetilde T}$, its label is the absolute value of the difference between its father's and grandfather's labels (see Figure~\ref{treeT}). 
Notice that all labels in $\mathbf{\widetilde T}$ are nonnegative. 

The random paths in the trees are coded by a sequence $(X_n)_{n\ge 3}$ of i.i.d. random variables taking values in the alphabet $\{R, L\}$ with probability $(p, 1-p)$. 
The path starts from the root and goes through its only child. Then the following steps are given by $(X_n)_{n\ge 3}$: 
Each $R$ corresponds to going through the right child (right step) and each $L$ corresponds to going through the left child (left step). 
Note that $F_n$ (respectively $\widetilde F_n$) is the label read in $\mathbf{T}$ (respectively $\mathbf{\widetilde T}$) at the end of the path $X_3\cdots X_n$. 

Observe that in the linear case, each step of the path can be interpreted as the right product of $(F_n, F_{n+1})$ by one of the two matrices 
\begin{equation}
\label{matrices}
A \egdef 
\begin{pmatrix}
0 & 1\\
1 & 1
\end{pmatrix}
\qquad \mbox{or}\qquad 
B \egdef 
\begin{pmatrix}
0 & -1\\
1 & 1
\end{pmatrix}.
\end{equation}
The non-linear case involves multiplication by matrices $A, B$ and 
$C\egdef 
\begin{pmatrix}
0 & 1\\
1 & -1
\end{pmatrix}$ but their distributions is no longer i.i.d., which makes this interpretation less convenient.

\subsection{Reduction of paths}

\begin{figure}[h]
	\begin{center}
	\input{reduction.pstex_t}
	\end{center}
\caption{The reduction pattern in the tree $\mathbf{T}$ (left) and in the tree $\mathbf{\widetilde T}$ (right).}
\label{reduction}
\end{figure}
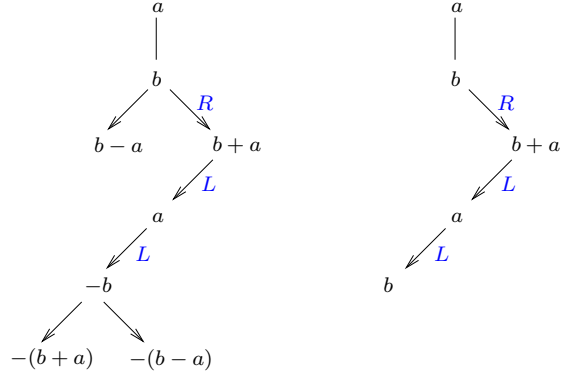

Our method relies on some properties of the trees $\mathbf{T}$ and $\mathbf{\widetilde T}$ illustrated in Figure~\ref{reduction}. 
In the following, we say an edge connecting a vertex $v$ and its child $v'$ is labelled by $(a, b)$ when $a$ is the label of $v$ and $b$ is the label of $v'$.

Suppose a random path goes through an edge labelled $(a, b)$ and then follows the pattern $RLL$ in $\mathbf{\widetilde T}$. Then it ends with an edge which is also labelled $(a, b)$. 
Therefore, we can remove from $(X_n)_n$ all occurences of the pattern $RLL$ when studying the sequence $(\widetilde F_n)_n$, as long as we keep in mind the number of such deletions. 

The linear case is a bit more complicated. 
Suppose now a random path in $\mathbf{T}$ goes through an edge labelled $(a, b)$. 
Notice the labels of the left and the right child are respectively $b-a$ and $b+a$. 
If the path follows the pattern $RLL$, then it ends on a vertex labelled by $-b$ and whose left child's and right child's labels are respectively $-(b+a)$ and $-(b-a)$. 
Since we are only interested in the behaviour of $|F_n|$, this allows us to remove in $(X_n)_n$ each pattern $RLL$, provided we exchange the following letter and keep in mind the number of deletions. 

This reduction process observed in $\mathbf{T}$ can be translated in the language of matrices by the following relations satisfied by $A$ and $B$:
\begin{eqnarray*}
ABBB & = & -A \label{RLLL}\\
ABBA & = & -B \label{RLLR}\, .
\end{eqnarray*}

\subsection{Reduced random Fibonacci sequence in the linear case.}
To formalize the reduction process, we associate to $X_3\ldots X_n$ a (generally) shorter word $W_n = Y_3^n\cdots Y_{k(n)}^n$ by the following induction. 
\begin{itemize}
\item $k(3) = 3$ and $Y_3^3 = X_3$. 
\item $W_{n+1}$ is deduced from $W_n$ in two steps. 

{\bf Step 1}: Add one letter ($R$ or $L$, see below) to the end of $W_n$. 

{\bf Step 2}: If the new word ends with the suffix $RLL$, remove this suffix. 

\noindent
Thus, we have either $W_{n+1}=W_n Y_{k(n)+1}^{n+1}$ (and $k(n+1)=k(n)+1$), or 
$W_{n+1} = Y_3^n\cdots Y_{k(n)-2}^n$ (and $k(n+1) = k(n)-2$). 
\end{itemize}
The letter which is added in step 1 depends on what happened when constructing $W_n$: 
\begin{itemize}
\item If $W_n$ was simply obtained by appending one letter to $W_{n-1}$ (or if $n=1$), we add $X_{n+1}$ to the end of $W_n$. 
\item Otherwise, we had removed the suffix $RLL$ when constructing $W_n$; we then add $\overline{X_{n+1}}$ to the end of $W_n$, where $\overline{R}\egdef L$ and $\overline{L}\egdef R$. 
\end{itemize}
Example: Let $X_3\dots X_9 = RRLLRLR$. Then, the successive reduced sequence $W_3, \dots, W_9$ are given by $R$, $RR$, $RRL$, $R$, $RL$, $\emptyset$, $L$. 

Observe that the label read in the tree $\mathbf{T}$ at the end of the path coded by $W_n$ has the same absolute value as $F_n$. 

\begin{lemma}
\label{Survival}
We denote by $|W_n|_R$ the number of $R$'s in $W_n$. We have
\begin{equation}
\label{RnR}
|W_n|_R\tend{n}{\infty}+\infty\qquad\mbox{a.s.}
\end{equation}
In particular, the length $k(n)$ of $W_n$ satisfies
$$
k(n) \tend{n}{\infty} +\infty.\qquad\mbox{a.s.}
$$
\end{lemma}

We postpone the proof of Lemma~\ref{Survival} to the end of the section. We will need in the sequel the following definition:
\begin{definition}
We say that $Y_{k(n)}^n$ \emph{survives} if, for all $m \ge n$, $k(m) \ge k(n)$. 
\end{definition}
The divergence of $k(n)$ in the previous lemma shows that, almost surely, for any $k\ge 3$, there exists a smallest $n_k$ such that $k(n_k)=k$ and $Y_{k}^{n_k}$ survives. 
In this case, $Y_{k}^{n_k} = Y_{k}^m$ for all $m\ge n_k$. We will then set 
$$ Y_k \egdef Y_{k}^{n_k}. $$
Note that $Y_{3} Y_4 \ldots$ contains infinitely many $R$'s (Lemma~\ref{Survival}), and no pattern $RLL$. Therefore the only place where consecutive $L$'s can appear is at the beginning of the sequence. However, these starting $L$'s are not relevant. 
If $Y_{3} Y_4 \ldots$ starts with an $L$, we can delete the first three letters without changing the values of the labels read along the path (see Figure~\ref{treeT}). Without loss of generality, we can henceforth assume that $Y_3=R$.

\begin{proof}[Proof of Lemma~\ref{Survival}]
We consider the successive changes in the number of $R$'s in $W_n$, which we denote by $S_1,S_2,\ldots\in\{\pm1\}$: If the $j$-th change corresponds to appending some $R$ to $W_n$, $S_j=1$; and if the $j$-th change is the deletion of some suffix $RLL$, $S_j=-1$ . Observe that if $S_j=1$, the only way to get $S_{j+1}=-1$ is to draw two $L$'s for the following two $X_n$'s (to remove the suffix $RLL$ before another $R$ is appended to $W_n$). Therefore,
$$
\PP (S_{j+1}=-1|S_j=1,S_{j-1},\ldots,S_1) = (1-p)^2. 
$$
On the other hand, $S_j=-1$ means that we have just deleted some suffix $RLL$, so that if the next $X_n$ is $L$ (which happens with probability $(1-p)$), we will have $S_{j+1}=1$. We thus get 
$$
\PP (S_{j+1}=1|S_j=-1,S_{j-1},\ldots,S_1) \ge (1-p). $$
We claim that the above requirements on the conditional distribution of $S_{j+1}$ imply
\begin{equation}
\label{1sur2-p}
\liminf_{L\to\infty}\dfrac{1}{L}\sum_{j=1}^L \ind{S_j=1} \ge \dfrac{1}{2-p} > \dfrac{1}{2}\quad\mbox{a.s.}
\end{equation}
This comes from the comparison between the stochastic process $(S_j)$ and the Markov chain $(M_\ell)$ taking values in $\{\pm1\}$, satisfying $M_1=1$ and having the following transition probabilities:
\begin{eqnarray*}
\PP(M_{l+1}=-1|M_\ell=1) & = & (1-p)^2,\\
\PP(M_{l+1}=1|M_\ell=-1) & = & (1-p).
\end{eqnarray*}
Note that the invariant probability measure for the Markov chain $(M_\ell)$ assigns mass  $\frac{1}{2-p}$ to the point $1$, so that
$$\lim_{L\to\infty}\dfrac{1}{L}\sum_{\ell=1}^L \ind{M_\ell=1} = \dfrac{1}{2-p} \mbox{ a.s.}$$

We are now going to extract a subsequence $(M_{\ell_j})$ from $(M_\ell)$ by deleting only some $-1$'s, and such that $(M_{\ell_j})$ has the same distribution as $(S_j)$. We set $\ell_1\egdef 1$. Suppose that $\ell_j$ is known. If $M_{\ell_j}=1$ or $M_{\ell_{j}+1}=1$, we set $\ell_{j+1}\egdef \ell_j+1$.  
Otherwise, $M_{\ell_j}=M_{\ell_{j}+1}=-1$. Let us denote by $\beta$ the probability $\PP (S_{j+1}=1|S_j=M_{\ell_j},S_{j-1}=M_{\ell_{j-1}},\ldots,S_1=M_{\ell_1})$. We know that $\beta\ge 1-p$. 
We then set
$$
\ell_{j+1}\ \egdef\ \begin{cases}
			\ell_j+1 & \mbox{with probability }\dfrac{1-\beta}{p},\\
			\inf\{k>\ell_j:M_k=1\} &  \mbox{with probability }\dfrac{\beta-1+p}{p}.
                 \end{cases}
$$
In this way, we get 
$$ \PP(M_{\ell_{j+1}}=1|M_{\ell_{j}}=-1,M_{\ell_{j-1}},\ldots,M_{\ell_1}) = \beta. $$
Since the proportion of $1$'s is greater in $(M_{\ell_j})$ than in $(M_\ell)$, \eqref{1sur2-p} follows, which in turn implies~\eqref{RnR}.
\end{proof}

\subsection{Reduced random Fibonacci sequence in the non-linear case.}
We associate to $X_3\ldots X_n$ the word $\widetilde W_n$, which is obtained by the same induction as $W_n$, except that the letter added in Step 1 is always $X_{n+1}$. 
The label read in the tree $\mathbf{\widetilde T}$ at the end of the path coded by $\widetilde W_n$ is equal to $\widetilde F_n$.

\begin{lemma}
\label{Survival2}
We denote by $|\widetilde W_n|_R$ the number of $R$'s in $\widetilde W_n$. 
If $p >1/3$, we have
$$
|\widetilde W_n|_R\tend{n}{\infty}+\infty\qquad\mbox{a.s.}
$$
In particular, the length of $\widetilde W_n$ goes to infinity almost surely.
\end{lemma}

\begin{proof}
Since each deletion of an $R$ goes with the deletion of two $L$'s, if $p>1/3$, the law of large numbers ensures that the number of remaining $R$'s goes to infinity. 
\end{proof}
The non-linear case for $p\le 1/3$ will be treated later (see Section~\ref{suite et fin}).

\subsection{Survival probability of an $R$}

We are now able to study both cases by introducing the probability $c$ of appending an $R$ after a deletion of the pattern $RLL$: $c=1-p$ in the linear case and $c=p$ in the non-linear case. In the sequel, we consider both the linear case for any $p\in ]0, 1]$ and the non-linear case for $p>1/3$. 
For simplicity, we will use the same notations $W_n=Y_3^n\dots Y_{k(n)}^n$ for the reduced word and $(Y_k)$ for the sequence of surviving letters in both cases. 

Observe that, by construction of the sequence $W_n$, if $Y_{k(n)}^n$ has been appended at time $n$, its survival only depends on the value of $Y_{k(n)}^n$ itself and the future $X_{n+1},X_{n+2}\ldots$.
We define 
$$ p_R\egdef\PP\left(Y_{k(n)}^n \mbox{ survives }|Y_{k(n)}^n=R \mbox{ has been appended at time }n\right).$$
A consequence of Lemma~\ref{Survival} and Lemma~\ref{Survival2} is that $p_R >0$. We now want to compute $p_R$ as a function of $p$. 

 We first need to analyze the following situation: Assume that in the construction of some $W_{n_0}$ we have deleted a suffix $RLL$. Then the survival of $Y_{k(n_0)}^{n_0}$ depends on the nature of what we call \emph{the next touching letter}, defined by the following algorithm:

{\bf Step 1} Set $n\egdef n_0$.

{\bf Step 2} If the letter appended to $W_n$ is $L$, it may interact with $Y_{k(n)}^n$ to delete it. Return $L$ as the next touching letter and halt the algorithm.

{\bf Step 3} If the letter appended to $W_n$ is $R$ and if this $R$ survives, then so does $Y_{k(n)}^n$, and we return $R$ as the next touching letter and halt the algorithm.

{\bf Step 4} Else, the letter appended to $W_n$ is a non-surviving $R$. Then there exists a smallest integer $m>n$ such that $k(m)=k(n)$ (corresponding to the time when this $R$ is deleted). Then set $n\egdef m$ and go back to Step~2.

Each time the algorithm enters Step~2, it has a probability $1-c$ to directly return $L$.
Since $1-c >0$, the algorithm ultimately halts with probability~1. 

We now go back to our computation of $p_R$.
Here is an exhaustive list of all cases in which $Y_{k(n)}^n=R$, appended at time $n$, survives:
\begin{itemize}
\item Case 1: $\mathbf{R}$ [\sout{$R\ldots$}] $\mathbf{R}$\\
$Y_{k(n+1)}^{n+1} = X_{n+1} = R$. Either it survives or it does not survive but the first touching letter after $\ell\ge 1$ deletions is an $R$; This happens with probability 
\begin{equation}
\label{p1}
p_{1} \egdef p p_R + p \sum_{\ell \ge 1} (1-p_R)^\ell c^\ell p_R = \frac{pp_R}{1-c(1-p_R)}\cdot
\end{equation}
\item Case 2: $\mathbf{R}$ [\sout{$R\ldots$}] $\mathbf{L}$ [\sout{$R\ldots$}] $\mathbf{R}$\\
Either $Y_{k(n+1)}^{n+1} = X_{n+1} = L$, or $Y_{k(n+1)}^{n+1} = X_{n+1} = R$ which does not survive but the first touching letter after $\ell\ge 1$ deletions is an $L$. In both cases, this $L$ is immediately followed either by a surviving $R$, or by a non-surviving $R$ and the second touching letter is an $R$; In view of the probability computed in Case 1, this happens with probability 
$$
p_{2} \egdef 
\left( (1-p) + p\sum_{\ell\ge 1} (1-p_R)^\ell c^{\ell - 1}(1-c)\right) \frac{p p_R}{1-c(1-p_R)} \cdot
$$
\end{itemize}
Writing $p_R = p_1 + p_2$, we get that $p_R$ satisfies
$$
p_R \Bigl( c^2 p_R^2 + (p^2+2c(1-p-c)) p_R +(1-c-2p)(1-c) \Bigr) = 0. 
$$
When $p_R\neq 0$, which is true in the linear case for any $p$ and in the non-linear case for $p>1/3$, this equation has only one non-negative solution given by 
\begin{equation}
\label{pR}
p_R = \frac{-p^2-2c(1-p-c)+p\sqrt{p^2+4c(1-p)}}{2c^2}.
\end{equation}
Observe that in the non-linear case for $p\le 1/3$, the above expression is non-positive. Thus, $p_R = 0$. 


\subsection{Distribution of the reduced sequence}
\label{distributionY}
We deduce from the preceding analysis the probability distribution, when $p_R>0$, of the reduced sequence $Y$ which is the concatenation of all surviving letters. 
By construction, since we assumed $Y_3=R$, consecutive $L$'s are not allowed in $Y$, hence
$$
\PP ( Y_{k+1} = R | Y_k = L) = 1.
$$
Moreover, if $Y_{k(n)}^n=R$ has been appended at time $n$ and survives, the following letter $Y_{k(n)+1}$ in $Y$ only depends on $X_{n+1}, X_{n+2}, \ldots$, which implies that $(Y_k)$ is a Markov chain. 
Observe that if $Y_k=R$, it is followed in $Y$ by another $R$ only in the last two cases of our exhaustive list. Hence, 
\begin{equation}
\label{defalpha}
\alpha \egdef \PP ( Y_{k+1} = R | Y_k = R) = \frac{p_1}{p_R} =
\frac{2c}{2c-p+\sqrt{p^2+4c(1-p)}}, 
\end{equation}
which leads to 
\begin{equation}
\label{alpha_both}
\alpha = 
\begin{cases}
\dfrac{3p-2+\sqrt{5p^2-8p+4}}{2p} & \mbox{ (linear case),}\\
\dfrac{2p}{p + \sqrt{p(4-3p)}}    & \mbox{ (non-linear case).}
\end{cases}
\end{equation}
The invariant probability measure of this Markov chain is given by 
\begin{equation}
\label{invariant_measure}
(\mu_R, \mu_L) \egdef \left( \frac{1}{2-\alpha}, \frac{1-\alpha}{2-\alpha} \right).
\end{equation}
>From now on, we denote by $\PP_\alpha$ the probability distribution on $\{R, L\}^{\{3, 4, \ldots\}}$ under which $(Y_k)$ is a Markov chain with $Y_3=R$,
$$
\PP_\alpha ( Y_{k+1} = R | Y_k = R)=\alpha \mbox{ and } \PP_\alpha ( Y_{k+1} = R | Y_k = L)=1\, .
$$

\subsection{Compression rate}
We are also interested in the ratio $\frac{k}{n_k}$ of surviving letters when $k$ is a large integer. 
The number $s_k$ of $R$'s in $Y_3\ldots Y_k$, which is the number of surviving $R$'s up to time $n_k$, satisfies $s_k/k = \mu_R +o(1)$. Let $d_k$ be the number of deleted $R$'s up to time $n_k$. 
Observe that the total number of $R$'s drawn up to time $n_k$ is $s_k + d_k$, so that $s_k/(s_k+d_k) = p_R +o(1)$. 
Since each deletion of an $R$ comes with the deletion of two $L$'s, we have $n_k = k+3d_k$. 
Therefore, 
\begin{equation}
\label{defsigma}
\frac{k}{n_k}\tend{k}{\infty}\sigma \egdef \left( 1 + 3\mu_R \frac{1-p_R}{p_R} \right)^{-1} \quad
\mbox{a.s}.
\end{equation}
It will be useful to see $\sigma$ as a function of $\alpha$. 
>From \eqref{p1}, we get $\alpha = p_1/p_R = \dfrac{p}{1-c(1-p_R)}$. Hence, 
$$
1-p_R = \frac{1}{c}\left(1-\frac{p}{\alpha}\right). 
$$
Moreover, \eqref{defalpha} yields 
$$
p = \alpha -c \frac{(1-\alpha)^2}{\alpha} \cdot
$$
Taking \eqref{invariant_measure} into account, an elementary computation leads to 
\begin{equation}
\label{alphatosigma}
\sigma= \frac{(2\alpha - 1)(2-\alpha)}{\alpha^2-\alpha+1} \cdot
\end{equation}

\section{Continued fractions in the tree $\mathbf{R}$}

\subsection{The tree $\mathbf{R}$}
The reduction of the sequence $(X_n)$ lead us, both in the linear and the non-linear case, to the study of a Markov chain $(Y_k)$ whose distribution is $\PP_\alpha$ on 
$\{R, L\}^{\{3, 4, \ldots\}}$. 
The only difference between the linear and the non-linear case is the value of the parameter $\alpha$. 

Let us consider the random sequence of integers $(G_k)_{k\ge 1}$, where $G_1=G_2=1$ and 
$G_{k}$ is the label read in the tree $\mathbf{T}$ ($\mathbf{\widetilde T}$ in the non-linear case) when following the path coded by $Y_3\cdots Y_k$. 
We get that for all $k\ge1$, $G_{k}$ has the same absolute value as $F_{n_k}$ ($\widetilde F_{n_k}$ in the non-linear case). 
We are thus left with the estimation of the exponential growth of the \emph{reduced Fibonacci sequence} $(G_k)_{k\ge 1}$. 

\begin{figure}
	\begin{center}
	\includegraphics[bb=0 0 124 208]{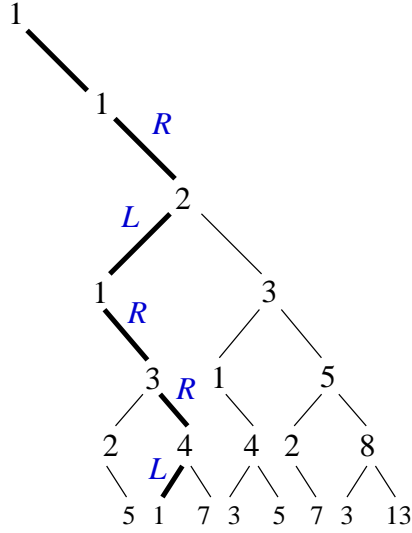}
	\end{center}
\caption{First lines of the tree $\mathbf{R}$. The random path marked with the bold edges is coded by the sequence $Y=RLRRL\ldots$ }
\label{treeR}
\end{figure}

Since $(Y_k)_k$ has no pattern $RLL$, the reduced Fibonacci sequence $(G_k)_{k\ge1}$ can be read along a random path in the tree $\mathbf{R}$ introduced in \cite{rittaud2006}, which is a sub-tree of $\mathbf{T}$ and $\mathbf{\widetilde T}$. 
The tree $\mathbf{R}$ is defined as follows: The root of $R$ has only one right child, which itself has only one right child. Any other vertex $v$ has either one right child or two (left and right) children, depending on whether $v$ is itself a left child or not: A left child has only one child, whereas a right child has two children. 
Each vertex is labelled with an integer : The root and its only child are labelled by $1$. 
If $v$ is a right child, its label is the sum of its father's and grandfather's labels; If $v$ is a left child, its label is the difference between its father's and grandfather's labels. 
(See Figure \ref{treeR}.) 

Please note that the step from the root to its child does not appear in $Y$: The condition $Y_3=R$ corresponds to the fact that the only child of the root has only one right child. For technical reasons, we will sometimes need to add an extra $R$ at the beginning of $Y$, representing the step from the root to its child, but this will be done explicitely.

\begin{figure}
	\begin{center}
	\input{probaR.pstex_t}
	\end{center}
\caption{Distribution probability of the random path in $\mathbf{R}$.}
\label{probaR}
\end{figure}
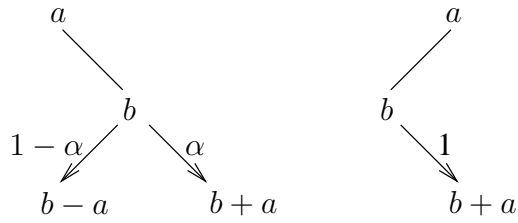

{From} the preceding section, we know that the distribution of the random path is given by a Markov chain: Each left step is followed by a right step, and each right step is followed by a right step with probability $\alpha$ and by a left step with probability $1-\alpha$, where $\alpha$ is given by~\eqref{defalpha} 
(see Figure~\ref{probaR}).

Let us recall some important properties of the tree $\mathbf{R}$. First, it is easily proved by induction that the label of each vertex is a \emph{positive} integer. As a consequence, we get that
$$\forall k\ge1,\quad G_k>0.
$$
Another easy induction shows that, if $x$ and $y$ are the labels of a child and its father, then $x$ and $y$ are relatively prime. Moreover, for each pair $(x, y)$ of relatively prime integers, there is exactly one vertex in $\mathbf{R}$ which is labelled by $x$ and its father by $y$ (see \cite{rittaud2006}). Therefore, for any positive rational number $q$, there is a unique vertex $x$ in $\mathbf{R}$ such that $q$ is the quotient of the label of $x$ by its father's. This vertex $x$ is a right child if and only if $q\ge1$.

\bigskip
For all $k\ge 3$, let $Q_k\egdef G_{k}/G_{k-1}$, so that 
\begin{equation}
\label{logGk}
\frac{1}{k} \log G_k = \frac{1}{k} \sum_{i=3}^{k} \log Q_i .
\end{equation}
The exponential growth of $G_k$ will thus be deduced from the probability distribution of $(Q_k)_{k\ge 3}$, which is related to the development in continued fractions of $Q_k$. 

\subsection{Shape of a path and continued fractions}

Let $\P$ be the set of all finite sequences $y=y_3 y_4\ldots y_k$ ($k\ge 3$) of $R$'s and $L$'s with $y_3=R$ and no pattern $LL$. 
Each sequence in $\P$ can be interpreted (with the same conventions as above) as a finite-length path in the tree $\mathbf{R}$.

To each sequence $y=y_3 y_4\ldots y_k \in\P$, we associate a rational number $q(y)$ defined as follows. 
We decompose the path $y_3 y_4\ldots y_k$ into pieces which are either elbows $RL$ or single right steps $R$. (This can be done in a unique way for all $y\in\P$.) Next, we introduce a cutting between each pair of (successive) identical pieces. We thus obtain a partition of the path into $\ell$ blocks; Let $a_1$ be the number of pieces in the last block, $a_2$ the number of pieces in the last but one block, and so on. 
If the last piece of the last block is an elbow, $q(y)$ is given by its development in continued fractions 
$$
q(y) \egdef [0, a_1,\ldots, a_\ell] = 0 + \cfrac{1}{a_1 + \cfrac{1}{a_2 + \cfrac{1}{\ddots+\cfrac{1}{a_\ell}}}}\ <1, 
$$ 
otherwise, $q(y)$ is set to 
$$
q(y) \egdef [a_1,\ldots, a_\ell] = a_1 + \cfrac{1}{a_2 + \cfrac{1}{\ddots+\cfrac{1}{a_\ell}}}\ >1. 
$$

The following proposition shows that the random variables $(Q_k)$ are precisely given by the preceding computation. 

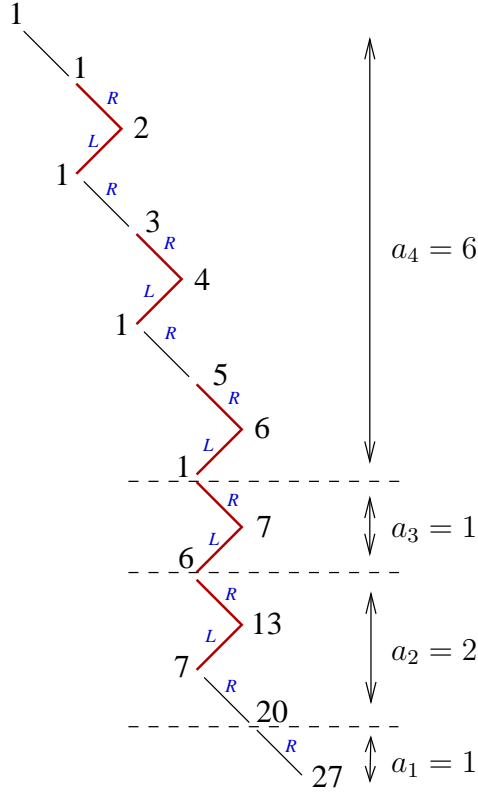
\begin{figure}[h]
	\begin{center}
	\input{ex_fraction.pstex_t}
	\end{center}
\caption{Consider the sequence $y=RLRRLRRLRLRLRR\in\P$. The corresponding path in the tree $\mathbf{R}$ is represented above, as well as $(g_j)_{1\le j\le 16}$. After decomposing $Ry$ into elbows and single right steps, and introducing the cuttings between identical pieces, we obtain the development in continued fractions of $q(Ry)=g_{16}/g_{15}=27/20$ by counting the number of pieces in each block: $27/20 = [1,2,1,6]$.}
\label{ex_fraction}
\end{figure}

\begin{prop}
Let $y=y_3\ldots y_k\in\P$, coding a finite-length path in $R$. 
Let $g_1=g_2=1,g_3, \ldots, g_{k}$ be the labels of the vertices visited by the path. 
Then $q(Ry) = g_{k}/g_{k-1}$. 
\end{prop}
(See Figure \ref{ex_fraction}.)

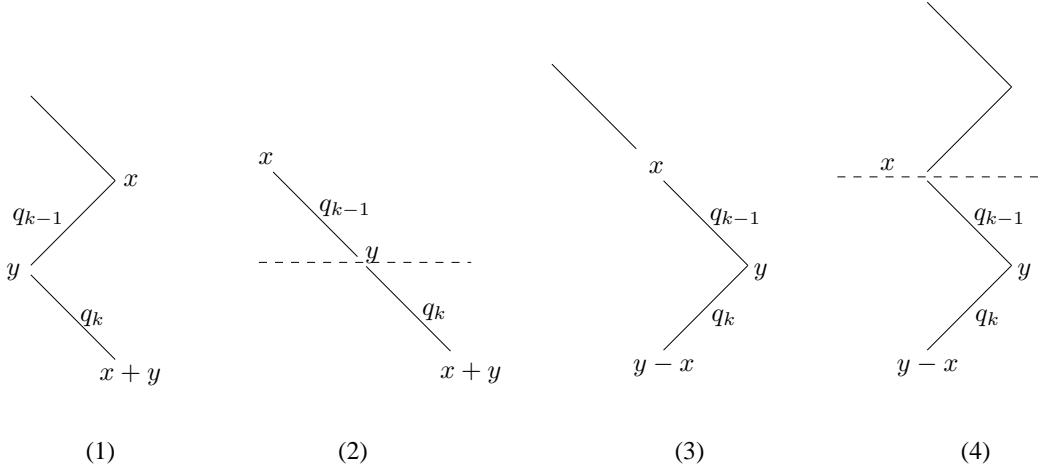
\begin{figure}
	\begin{center}
	\input{evolvQ.pstex_t}
	\end{center}
\caption{From $q_{k-1}$ to $q_{k}$. (1) corresponds to the case when $q_{k-1}<1$, (2), (3) and (4) correspond to the cases when $q_{k-1} >1$.}
\label{Q}
\end{figure}

\begin{proof}
We proceed by induction on $k$. Let us denote by $q_k$ the quotient $g_k/g_{k-1}$. 
For $k=3$, the sequence $Ry$ is reduced to $RR$: We thus have two single right steps separated by a cutting and we check that $$q_3 = \frac{2}{1} = 1 + \cfrac{1}{1}=q(RR).$$ 
Assume the result is proved up to $k-1$ and consider a sequence $y$ of $R$ and $L$'s of length $k$. Let us analyze the three possible configurations for $y_{k-1}y_{k}$. 
\begin{itemize}
\item $y_{k-1}y_{k}=LR$ (see Figure \ref{Q}, case (1)): Since the last letter is an $R$, we have 
$$
q_{k} = \frac{g_{k}}{g_{k-1}} = \frac{g_{k-1}+g_{k-2}}{g_{k-1}} = 1 + \frac{1}{q_{k-1}}\ .
$$
On the other hand, the induction hypothesis gives 
$q_{k-1} = q(Ry_3\ldots y_{k-1})=[0, a_1, a_2, \ldots, a_\ell]$, with $(a_i)$ given by the number of pieces in the blocks. Hence, 
$q_{k} = [a_1+1, a_2, \ldots, a_\ell]$. 
Since there is no change in the cuttings when appending the last $R$, the number $a_1$ of pieces in the last block is increased by 1 and the other ones are left unchanged. Therefore, 
$q(Ry_3\ldots y_{k})=q_{k}$. 
\item $y_{k-1}y_{k}=RR$ (see Figure \ref{Q}, case (2)): We still have $q_{k} = 1 + 1/q_{k-1}$. The difference with the previous case is that $q_{k-1} =[a_1, a_2, \ldots, a_\ell] > 1$, thus $q_{k} = [1, a_1, a_2, \ldots, a_\ell]$. 
Here, the last $R$ has introduced a new cutting, so that we have one more block of length 1 and the other blocks are left unchanged. Therefore, 
$q(Ry_3\ldots y_{k})=q_{k}$. 
\item $y_{k-1}y_{k}=RL$ (see Figure \ref{Q}, cases (3) and (4)): Since the last letter is an $L$, we have 
$$
q_{k} = \frac{g_{k}}{g_{k-1}} = \frac{g_{k-1}-g_{k-2}}{g_{k-1}} = 
\cfrac{1}{1+\cfrac{1}{q_{k-1} - 1}}\ .
$$
The induction hypothesis gives that $q_{k-1} = q(Ry_3\ldots y_{k-1})=[a_1, a_2, \ldots, a_\ell]$, with $(a_i)$ given by the number of pieces in the blocks. Hence, the development in continued fractions of $q_{k}$ depends on the value of $a_1$. 

If $a_1=1$ (case (3)), we get $q_{k} = [0, a_2+1, a_3, \ldots, a_\ell]$; 
Appending the last $L$ transforms the last piece into an elbow. The fact that $a_1=1$ means we had a cutting just before this last piece, which disappears after the transformation. Therefore, we have one less block and the number of pieces of the last remaining block is increased by 1. 

If $a_1> 1$ (case (4)), we get $q_{k} = [0, 1, a_1-1, a_2, \ldots, a_\ell]$; The fact that $a_1>1$ means we had no cutting just before the last piece, and one is created when appending the last $L$. Therefore, we have one more block of length 1, and the number of pieces of the last but one block is decreased by 1. 

In both cases, we get $q(Ry_3\ldots y_{k})=q_{k}$.
\end{itemize}
\end{proof}

\subsection{Stern-Brocot intervals}
We recall here the definition of the positive \emph{Stern-Brocot intervals}, which are the real intervals of the form $[\frac{a}{b},\frac{c}{d}]$ where $a,b,c,d$ are nonnegative integers with \hbox{$ad-bc=-1$}. (We adopt the usual convention that $\frac{1}{0}=+\infty$.) All these intervals can be obtained inductively, starting with the interval $[\frac{0}{1},\frac{1}{0}]$ (the only Stern-Brocot interval of rank 0). Each Stern-Brocot interval $[\frac{a}{b},\frac{c}{d}]$ of rank $r$ is cut into two Stern-Brocot intervals of rank $(r+1)$: $[\frac{a}{b},\frac{a+c}{b+d}]$ and $[\frac{a+c}{b+d},\frac{c}{d}]$. Recall also that $\frac{a+c}{b+d}$ is called the \emph{mediant} of $\frac{a}{b}$ and $\frac{c}{d}$.

\medskip
Let us fix a sequence $s=s_1\dots s_{m} \in \P$ (we also include here the case where $s$ is the empty sequence).
We now consider the set of finite paths in the tree $\mathbf{R}$ having $s$ as a suffix:
$$
\Y_{s} \egdef 
\left\{ Ry: y=(y_i)_{3\le i \le k} \in\P,\ k\ge m+2 ;\ y_{k-m+1}\ldots y_k = s_1\dots s_{m} \right\}\ .
$$

\begin{prop}
\label{SB}
Let $r$ be the number of pieces in the decomposition of $s$ into elbows $RL$ and single right steps $R$. Then
\begin{itemize}
\item The closure $\overline{q(\Y_s)}$ of $q(\Y_s)$ is a Stern-Brocot interval of rank $r$;
\item if $r$ is even (respectively odd), $\overline{q(\Y_{RLs})}$ is the left (respectively right) Stern-Brocot sub-interval of rank $r+1$ of $\overline{q(\Y_s)}$, and $\overline{q(\Y_{Rs})}$ is the right (respectively left) one.
\end{itemize}
\end{prop}

\begin{proof}
We proceed by induction on $r$. If $r=0$, $s$ is the empty sequence and $\Y_s$ is the set of all finite paths in $R$, so $\overline{q(\Y_s)}=[0,+\infty]$. Moreover, $\Y_{Rs}=\Y_{R}$ is the set of all finite paths ending with an $R$, hence $\overline{q(\Y_{Rs})}=[1,+\infty]$. $\Y_{RLs}=\Y_{RL}$ is the set of all finite paths ending with an elbow and $\overline{q(\Y_{RLs})}=[0,1]$. Therefore $\overline{q(\Y_{RLs})}$ is the left Stern-Brocot sub-interval of rank~1.

We consider a suffix $s$ such that $\overline{q(\Y_s)}$ is a Stern-Brocot interval of rank $r$, and its two possible extensions $Rs$ and $RLs$. 
Assume for simplicity that $s$ ends with an $R$ (the proof is the same when $s$ ends with $RL$). 
Let $\ell$ be the number of blocks in $s$: $q(s)=[a_1, \ldots, a_\ell]$.
Since the last letters of $y$ give the beginning of the development in continued fractions of $q(y)$, all $q(y)$ for $y\in \Y_{s}$ have their first $(\ell - 1)$ partial quotients fixed, and equal to those of $q(s)$. Moreover, their $\ell$-th partial quotient is at least $a_\ell$. 
Thus, $\overline{q(\Y_s)}$ is the interval whose bounds are $[a_1, \ldots, a_{\ell-1}]$ and $[a_1, \ldots, a_{\ell}]$. 
Notice the way these bounds are ordered depends on the parity of $\ell$: 
$[a_1, \ldots, a_{\ell-1}]<[a_1, \ldots, a_{\ell}]$ iff $\ell$ is even. 

Either $Rs$ (if the first piece of $s$ is an $R$) or $RLs$ (if the first piece of $s$ is an $RL$) counts one more cutting than $s$ and gives rise to the interval whose bounds are $[a_1, \ldots, a_{\ell -1}, a_\ell, 1]$ and $[a_1, \ldots, a_{\ell -1}, a_\ell]$. The other one gives rise to the interval whose bounds are $[a_1, \ldots, a_{\ell -1}]$ and $[a_1, \ldots, a_{\ell -1}, a_\ell + 1]$. In fact, since $[a_1, \ldots, a_{\ell -1}, a_\ell, 1]=[a_1, \ldots, a_{\ell -1}, a_\ell + 1]$, the two intervals have a common bound. We let the reader check that this common bound is the mediant of $[a_1, \ldots, a_{\ell-1}]$ and $[a_1, \ldots, a_{\ell}]$, so that $\overline{q(\Y_{Rs})}$ and $\overline{q(\Y_{RLs})}$ are two Stern-Brocot intervals of rank~$r+1$. 
As a consequence of the previous remarks, the table below gives the relative positions of $\overline{q(\Y_{RLs})}$ and $\overline{q(\Y_{Rs})}$:

\begin{center}
\begin{tabular}{|c|c|c|}
\hline 
 & the first piece of $s$ is $R$ &  the first piece of $s$ is $RL$ \\ 
\hline 
$\ell$ odd  &  $Rs\leftrightarrow$ left   & $RLs\leftrightarrow$ left \\
	    &  $RLs\leftrightarrow$ right & $Rs\leftrightarrow$  right \\ 
\hline 
$\ell$ even &  $RLs\leftrightarrow$ left  & $Rs\leftrightarrow$ left   \\
	    &  $Rs\leftrightarrow$ right  & $RLs\leftrightarrow$  right \\
\hline 
\end{tabular}
\end{center}
Each time we extend the suffix $s$ with one more piece, either this piece is similar to the first piece of $s$, which introduces a new cutting, and $\ell$ is increased by $1$, or the new piece is different and $\ell$ is unchanged. Therefore, we move either vertically or horizontally in the table, so that the relative positions of $\overline{q(\Y_{Rs})}$ and $\overline{q(\Y_{RLs})}$ alternate. 
\end{proof}

\subsection{Probability distribution of $(Q_k)_{k\ge 3}$}
\label{section_distributionQ}
We now turn back to the Markov chain $Y$ following the probability distribution $\PP_\alpha$ (see Section~\ref{distributionY}). 
The ergodic theorem for this Markov chain gives that, almost surely,
$$ 
\dfrac{1}{k}\sum_{i=3}^k\ind{RY_3\ldots Y_i\in\Y_R}\tend{k}{\infty} \mu_R.
$$
If we decompose the sequence $Y_3\ldots Y_i$ into pieces $RL$ and $R$, it is not hard to see that all but the last piece appear independently, with probability $1-\alpha$ for $RL$ and $\alpha$ for $R$. Therefore, if we fix some $s\in\P$ and denote by $|s|_{RL}$ (respectively $|s|_{R}$) the number of pieces $RL$ (respectively $R$) in the decomposition of $s$ into pieces, we get from the law of large numbers that
\begin{equation}
\label{lln}
\dfrac{1}{k}\sum_{i=3}^k\ind{RY_3\ldots Y_i\in\Y_s}\tend{k}{\infty} c(s)\,\alpha^{|s|_{R}}(1-\alpha)^{|s|_{RL}},
\end{equation}
where
$$
c(s) = \begin{cases}
	\mu_R/\alpha & \mbox{if $s$ ends with an $R$,}\\
	(1-\mu_R)/(1-\alpha) & \mbox{otherwise}.
       \end{cases}
$$

Since for all $i\ge 3$ we have $Q_i=q(RY_3\ldots Y_{i})$, it is natural to introduce the following probability distribution
$\nu_\alpha$ on $\RR_+$: $\nu_\alpha$ is defined by
$$
\forall s\in\P,\quad \nu_\alpha\left(\overline{q(\Y_{s})}\right)\egdef \alpha^{|s|_{R}}(1-\alpha)^{|s|_{RL}}.
$$
In view of Proposition~\ref{SB}, this amounts to define it inductively on Stern-Brocot intervals in the following way: First assign mass $1-\alpha$ to $[0,1]$ and $\alpha$ to $[1,\infty]$. Once $\nu_\alpha$ is defined on some Stern-Brocot interval of rank $r$, a proportion $\alpha$ of its mass is given to the left (respectively right) subinterval of rank $r+1$ when $r$ is odd (respectively even) (See Figure~\ref{mesure}). We can notice the similarity between this construction and the Denjoy-Minkowski measure $\mu^{(\alpha)}$ presented in \cite{chassaing1984}. The difference lies in the fact that the proportion $\alpha$ is in our case alternatively given to the left and the right subinterval.

{From} \eqref{lln}, we obtain that, for all $f\in L^1(\nu_\alpha)$,
\begin{equation}
\label{distributionQ}
\dfrac{1}{k}\sum_{i=3}^k f(Q_i) \tend{k}{\infty} 
\dfrac{1-\mu_R}{1-\alpha} \int_0^1 f(x)\, d\nu_\alpha(x) + 
\dfrac{\mu_R}{\alpha} \int_1^\infty f(x)\, d\nu_\alpha(x)\quad\mbox{a.s.}
\end{equation}

\begin{remark}
\label{remark nu_alpha}
Observe that we need the correction $c(s)$ in \eqref{lln} because the last piece of $Y_3\ldots
Y_i$ has a different distribution. If we first cut the infinite sequence $(Y_i)_{i\ge 3}$ into pieces $R$ and $RL$, and denote by $i_j$ the index at the end of the $j$-th piece, then the distribution of $q(RY_3\ldots Y_{i_j})$ converges to $\nu_\alpha$ without correction.
\end{remark}

\section{The Lyapunov exponent}

\subsection{Computation of the Lyapunov exponent}

We now end the proof of Theorem~\ref{MainTheorem}. 
It is easily seen that $x\mapsto\log x$ belongs to $L^1(\nu_\alpha)$. 
Using \eqref{logGk} and \eqref{distributionQ}, we obtain
\begin{equation*}
\dfrac{1}{k}\log G_k \tend{k}{\infty} 
\dfrac{1-\mu_R}{1-\alpha} \int_0^1 \log(x)\, d\nu_\alpha(x) + 
\dfrac{\mu_R}{\alpha} \int_1^\infty \log(x)\, d\nu_\alpha(x)
\quad\mbox{a.s.}
\end{equation*}
Using \eqref{defsigma} and $\mu_R=1/(2-\alpha)$, we get 
\begin{equation}
\label{presquefin}
\dfrac{1}{n_k} \log |F_{n_k}| \tend{k}{\infty} 
\dfrac{\sigma}{2-\alpha}  \left( \int_0^1 \log(x)\, d\nu_\alpha(x) + 
\dfrac{1}{\alpha} \int_1^\infty \log(x)\, d\nu_\alpha(x)\right)\quad\mbox{a.s.}
\end{equation}
In the linear case, since we are dealing with a product of i.i.d. matrices, we know that the limit
of $(1/n)\log |F_n|$ exists almost surely, and is given by the largest Lyapunov exponent.

\medskip

Of course, we get the same formula as \eqref{presquefin} for $(1/n_k) \log \widetilde F_{n_k}$ in the non-linear case for $p>1/3$ (where $\alpha$ is given by \eqref{alpha_both}). 
As we already pointed out, this case cannot be reduced to a product of i.i.d. matrices, therefore we need a little argument to get the almost-sure existence of the limit of $(1/n)\log \widetilde F_n$. It consists in controlling the size of the deleted blocks.

We already know the almost-sure convergence of $(1/n_k)\log \widetilde F_{n_k}$, along the subsequence $(n_k)$ corresponding to surviving letters after the reduction process. 
Consider now  $n$ lying between $n_k$ and $n_{k+1}$. We have
$$
\left|\dfrac{1}{n}\log \widetilde F_{n} - \dfrac{1}{n_k}\log \widetilde F_{n_k} \right|\le \dfrac{n-n_k}{n} \dfrac{1}{n_k}\log \widetilde F_{n_k} + \dfrac{1}{n}\left|\log \dfrac{\widetilde F_{n}}{\widetilde F_{n_k}}\right|.
$$
We need to control the quantity $T_k\egdef n_{k+1}-n_k$, which is 1 plus the number of deleted letters between two successive surviving letters. The probability distribution of the random variable $T_k$ is given by the law of the following stopping time for the i.i.d. sequence on the alphabet $\{R,L\}$ with probability $(p,1-p)$: 
Draw a sample of this process and stop the first time the number of $L$'s is equal to twice the number of $R$'s plus one. 
This sample without the last $L$ corresponds to all the patterns $RLL$, between two successive surviving letters, which are removed during the reduction process. 
When $p>1/3$, the stopping time is almost surely finite and its expected value is finite. 

Since the $T_k$'s are i.i.d. and have a finite expected value, we get
$$
\dfrac{n-n_k}{n} \le \dfrac{T_k}{k} \tend{k}{\infty} 0 \quad\mbox{a.s.}
$$
Observing that $\widetilde F_{n} \le \max\{\widetilde F_{n_k},\widetilde F_{n_{k-1}}\}2^{n-n_k}$, the convergence along the subsequence $(n_k)$ is enough to conclude that the limit of $(1/n)\log \widetilde F_{n}$ exists almost surely, and is given by the right-handside of~\eqref{presquefin}.

\medskip
It remains now to prove that this limit is equal to $\int_0^1 \log(x)\, d\nu_\alpha(x)$. To this end, we use some changes of variables in the computation of the integral.

\begin{lemma}
\label{changement_variable}
For all $f\in L^1(\nu_\alpha)$, 
\begin{eqnarray*}
\int_0^1 f\left(\frac{1}{1-x}\right) d\nu_\alpha(x) &=& \frac{1-\alpha}{\alpha}\int_1^\infty f(x) d\nu_\alpha(x)\\
\int_0^1 f\left(\frac{1-x}{x}\right) d\nu_\alpha(x) &=& (1-\alpha) \int_0^\infty f(x) d\nu_\alpha(x)
\end{eqnarray*}
\end{lemma}

\begin{proof}
It is sufficient to prove the equalities when $f=\ind{[a/b, c/d]}$ is the indicator function of a Stern-Brocot interval ($ad-bc=-1$).
The first equality becomes, for $a\ge b$, 
$$
\nu_\alpha\left(\left[\frac{a-b}{a}, \frac{c-d}{c}\right]\right) 
= \frac{1-\alpha}{\alpha}\ \nu_\alpha\left(\left[\frac{a}{b}, \frac{c}{d}\right]\right) .
$$
Observe that $\left[\frac{a-b}{a}, \frac{c-d}{c}\right]$ is also a Stern-Brocot interval, which has the same rank $r$ as $[\frac{a}{b}, \frac{c}{d}]$. 
For $r=1$, the equality holds trivially. 
Since $x\mapsto 1/(1-x)$ is increasing, the right sub-interval of $[\frac{a}{b}, \frac{c}{d}]$ corresponds to the right sub-interval of $\left[\frac{a-b}{a}, \frac{c-d}{c}\right]$. The result follows by induction on $r$. 

For $f=\ind{[\frac{a}{b}, \frac{c}{d}]}$, the second equality reduces to 
$$
\nu_\alpha\left(\left[\frac{d}{c+d}, \frac{b}{a+b}\right]\right) 
= (1-\alpha)\ \nu_\alpha\left(\left[\frac{a}{b}, \frac{c}{d}\right]\right) .
$$
If $[\frac{a}{b}, \frac{c}{d}]$ is a Stern-Brocot interval of rank $r$, 
$\left[\frac{d}{c+d}, \frac{b}{a+b}\right]$ is also a Stern-Brocot interval, but of rank $r+1$. 
For $r=1$, the equality holds trivially. 
Since $x\mapsto (1-x)/x$ is decreasing, the right sub-interval of $[\frac{a}{b}, \frac{c}{d}]$ corresponds to the left sub-interval of $\left[\frac{d}{c+d}, \frac{b}{a+b}\right]$. The result follows by induction on $r$. 
\end{proof}

Applying Lemma~\ref{changement_variable} to $x\mapsto\log x$ and summing, we get that 
$$
- \int_0^1 \log (x) d\nu_\alpha(x) = 
\frac{1-\alpha}{\alpha}\int_1^\infty \log(x) d\nu_\alpha(x) +
(1-\alpha) \int_0^\infty \log(x) d\nu_\alpha(x) . 
$$
It follows immediately that 
\begin{eqnarray}
\int_0^1 \log (x) d\nu_\alpha(x)&=& \frac{(\alpha + 1)(1-\alpha)}{1-2\alpha}\int_0^\infty \log(x) d\nu_\alpha(x) ,\nonumber\\
\int_1^\infty \log (x) d\nu_\alpha(x)&=& \frac{\alpha(\alpha-2)}{1-2\alpha}\int_0^\infty \log(x) d\nu_\alpha(x) \label{gamma_positif}.
\end{eqnarray}
Substituting in the right-handside of \eqref{presquefin} and recalling \eqref{alphatosigma}, we conclude the proof of Theorem~\ref{MainTheorem}, for all $p$ in the linear case and for $p>1/3$ in the non-linear case. The non-linear case for $p\le 1/3$ is treated in the next section. 

\subsection{Variation properties of the Lyapunov exponent}
\label{suite et fin}
We now want to prove that $p\mapsto \widetilde \gamma_p$ (non-linear case) is a non-decreasing function. We first establish a comparison lemma.

\begin{lemma}
\label{compare}
Let $x$ be a path in the tree $\mathbf{\widetilde T}$, and let $x'$ be obtained from $x$ by turning an $L$ into an $R$. Then any label read along $x$ is always smaller than the corresponding label read along $x'$.
\end{lemma}

\begin{proof}
Let $y=y_3\dots y_k\in\P$ coding the end of a finite-length path in $\mathbf{R}$. 
Assume that the vertices of the edge preceding $y$ are labelled by $a$ and $b$. 
Then the last vertex of the path is labelled by a linear combination of $a$ and $b$ of the form 
$d(y)a+n(y)b$, where $n(y)$ and $d(y)$ are integers depending only on $y$ (see Figure~\ref{croissanceR}, left).
More formally, $n(y)$ and $d(y)$ can be defined by the following induction. 
$n(\emptyset)\egdef 1$, $d(\emptyset)\egdef 0$, $n(R)\egdef d(R)\egdef 1$, 
and for $i\ge 3$
\begin{eqnarray*}
n(y_3\dots y_i R) & \egdef & n(y_3\dots y_i) + n(y_3\dots y_{i-1})\\
d(y_3\dots y_i R) & \egdef & d(y_3\dots y_i) + d(y_3\dots y_{i-1})\\
n(y_3\dots y_i L) & \egdef & n(y_3\dots y_i) - n(y_3\dots y_{i-1})\\
d(y_3\dots y_i L) & \egdef & d(y_3\dots y_i) - d(y_3\dots y_{i-1})\, .
\end{eqnarray*}
Since $y$ codes a path in $\mathbf{R}$, an induction shows that $d(y)$ and $n(y)$ are nonnegative. 
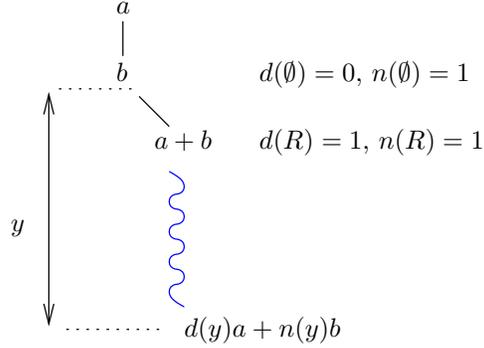
\begin{figure}[h]
\input{def_n_et_d.pstex_t}
\caption{Definition of $d(y)$ and $n(y)$ (left); Comparison between two paths in $\mathbf{R}$ (right).}
\label{croissanceR}
\end{figure} 


Consider two paths $x$ and $x'$ in $\mathbf{\widetilde T}$ differing at level $i$ : We decompose the end of the paths from level $i$ as $LL^ry$ and $RL^ry$, where $y$ starts with an $R$ and $r\ge0$. 

Suppose first that, after the difference, all letters of $x$ and $x'$ are $L$'s ($y=\emptyset$). 
We let the reader check that the labels after level $i$ and $i+1$ are well-ordered. 
Moreover, the label after level $k\ge i+2$ in $x'$ is equal to the label after level $k-3$ in $x$. This can be seen by making a reduction (removing a pattern $RLL$) in the path $x'$. 
Denoting by $a$ and $b$ the (nonnegative) labels after level $k-3$ and $k-2$ in $x$, the label after level $k$ in $x$ is given by $|b-|b-a||\le a$. 

Suppose now that the suffix $y$ is reduced ($y\in\P$). The above argument shows that labels after level $i+r-1$ and $i+r$ are well-ordered: Denote these labels by $a$ and $b$ in $x$, and $a'$ and $x'$ in $x'$. Then the label at the end of $x$ (respectively $x'$) is the linear combination $d(y)a + n(y)b$ (respectively $d(y)a' + n(y)b'$). Since $d(y)$ and $n(y)$ are nonnegative, we conclude this case.

In the general case, we make all possible reductions on $y$. We are left either with a sequence in $\P$ or with a sequence of $L$'s, which are the two situations we have already studied.

\end{proof}

\begin{prop}
If $p\le1/3$, $\widetilde \gamma_p\egdef \lim_n (1/n)\log \widetilde F_n$ exists and is equal to 0. The function $p\mapsto \widetilde \gamma_p$ is continuous and non-decreasing.
\end{prop}

\begin{proof}
We start by proving that $\lim_{p\searrow 1/3}\widetilde \gamma_p = 0$. 
We can extend the definition of $\nu_\alpha$ to any $\alpha\in [0,1]$, and also set for $0\le\alpha\le 1/2$
$$ \gamma(\alpha)\egdef\int_0^\infty \log(x) d\nu_\alpha(x). $$
Observe that, if $[a,b]$ is a Stern-Brocot interval, then $[1/b,1/a]$ is also a Stern-Brocot interval, satisfying
$$ \nu_\alpha\Bigl([a,b]\Bigr) = \nu_{1-\alpha}\Bigl([1/b,1/a]\Bigr). $$
Therefore, we get
$$ \gamma(\alpha)=-\gamma(1-\alpha). $$
In particular, $\gamma(1/2)=0$. Moreover, $\alpha\mapsto \gamma(\alpha)$ is easily seen to be a continuous function (we can write it as a uniform limit of continuous functions). Hence
$$ \lim_{p\searrow 1/3}\widetilde \gamma_p = \lim_{\alpha\to 1/2} \gamma(\alpha) = 0. $$

Now, let $0<p\le p'\le1$. Let $(X_n)$ and $(X'_n)$ be random paths in $\mathbf{\widetilde T}$ for the respective parameters $p$ and $p'$. We can realize a coupling of $(X_n)$ and $(X'_n)$ such that for any $n$, $X_n=R$ implies $X'_n=R$. From Lemma~\ref{compare}, it follows that the label $\widetilde F_n$ read along $X$ is always smaller than the label $\widetilde F'_n$ read along $X'$. If we choose $p\le 1/3$ and $p'>1/3$, we get that
$$\limsup \dfrac{1}{n}\log \widetilde F_n\le \lim\dfrac{1}{n}\log \widetilde F'_n = \widetilde\gamma_{p'}.$$
Since $\lim_{p'\searrow 1/3}\widetilde \gamma_{p'} =0$, we deduce that $\widetilde \gamma_{p} =0$ for any $p\le1/3$.
Moreover, this argument obviously shows that $p\mapsto \widetilde \gamma_p$ is a non-decreasing function.
\end{proof}

\begin{corollary}
The function $p\mapsto \gamma_p$ is increasing on $]0, 1[$. The function $p\mapsto \widetilde\gamma_p$ is increasing on $]1/3,1]$.
\end{corollary}

\begin{proof}
Recall that $\gamma(\alpha) = \int_0^\infty\log x d\nu_\alpha(x)=\widetilde\gamma_{p}$, for $\alpha= 2p/(p + \sqrt{p(4-3p)})$. The function $p\mapsto 2p/(p + \sqrt{p(4-3p)})$ is increasing and sends $]1/3,1]$ onto $]1/2,1]$. Hence $\gamma(\alpha)$ is non-decreasing on $]1/2,1]$. 

In the linear case, we also have $\gamma_{p}=\gamma(\alpha)$ where now $\alpha= (3p-2+\sqrt{5p^2-8p+4})/(2p)\in]1/2,1]$. Since this expression is also increasing in $p$, $p\mapsto \gamma_p$ is non-decreasing. 
Moreover, we easily deduce from \eqref{gamma_positif} that $\gamma_p>0$ for any $p\in]0, 1[$.
We also know from \cite{peres1990} that $\gamma_{p}$ is an analytic function of $p\in]0, 1[$, thus it is increasing.

This in turn implies that $\gamma(\alpha)$ is increasing on $]1/2,1]$, and so is $\widetilde\gamma_{p}$ for $p\in]1/3,1]$.
\end{proof}

\subsection{Derivative of the Lyapunov exponents}
\label{study}
%
%
%

The following proposition gives a formula for the derivative of $\gamma$ with respect to $\alpha$, which uses the product measure $\nu_{\alpha}\otimes\nu_{\alpha}$.
\begin{prop}
\label{prop_derivee_gamma}
For all $1/2< \alpha\le 1$, 
\begin{equation}
\label{derivee_gamma}
\gamma'(\alpha)
= \gamma(\alpha)\frac{1+2\alpha -2\alpha^2}{(2\alpha-1)(\alpha^2-\alpha+1)}
 + \frac{2\alpha-1}{\alpha^2-\alpha+1}\int_0^{\infty}\int_0^{\infty}\log\frac{x+y+xy}{x+y+1}\,  \nu_{\alpha}\otimes\nu_{\alpha}(dx,dy)\, .
\end{equation}
\end{prop}

Before proving this formula, we now use it to compute the derivatives of the Lyapunov exponents for $p=1$. 
%
When $p=1$, $\alpha=1$ (both in the linear and in the non-linear cases), and $\nu_\alpha$ is the Dirac measure on $\phi$. 
\eqref{derivee_gamma} yields 
$$
\frac{d\gamma}{d\alpha}(1) = \log\phi + \log\frac{2\phi+\phi^2}{2\phi+1} = \frac{\log 5}{2}\, .
$$
Since $\dfrac{d\alpha}{dp}(1)=1$ (both in the linear and in the non-linear cases), we get $$\dfrac{d\gamma_p}{dp}(1) = \dfrac{d\widetilde\gamma_p}{dp}(1) = (\log 5)/2. $$

We now turn to the proof of Proposition \ref{prop_derivee_gamma}.

\begin{proof}[Proof of Proposition~\ref{prop_derivee_gamma}]
We fix $1/2<\alpha-\varepsilon<\alpha\le1$. Let $Y$ be a Markov chain following $\PP_\alpha$, decomposed into pieces $R$ and $RL$. We decide independently to change each piece $R$ into $RL$ with probability $\varepsilon/\alpha$. We thus obtain a new Markov chain $Y'$ following $\PP_{\alpha-\varepsilon}$.
Let $(G_k)_{k\ge1}$ and $(G'_k)_{k\ge1}$ be the labels read along the paths $Y$ and $Y'$ respectively. We introduce the subsequences $(k_m)$ and $(k'_m)$ such that $G_{k_m}$ and $G'_{k'_m}$ are the labels read after the $m$-th piece of $Y$ and $Y'$ respectively. 
We know that 
\begin{eqnarray}
\label{on verra}
\gamma(\alpha) - \gamma(\alpha-\varepsilon)
&=& \sigma(\alpha)\lim_{m\to\infty}  \EE\left[\frac{1}{k_m}\log G_{k_m}\right]-\sigma(\alpha-\varepsilon)\lim_{m\to\infty}  \EE\left[\frac{1}{k'_m}\log G'_{k'_m}\right]\nonumber\\
&=& (\sigma(\alpha)-\sigma(\alpha-\varepsilon))\dfrac{\gamma(\alpha)}{\sigma(\alpha)}
+\sigma(\alpha-\varepsilon) \lim_{m\to\infty}\EE\left[\frac{1}{k_m}\log G_{k_m}-\frac{1}{k'_m}\log G'_{k'_m}\right].
\end{eqnarray}
Since
$$ \dfrac{k_m}{m}\tend{m}{\infty}2-\alpha\mbox{ and }\dfrac{k'_m}{m}\tend{m}{\infty}2-\alpha+\varepsilon\mbox{ a.s.},
$$
the last term in \eqref{on verra} becomes
$$ \frac{\sigma(\alpha-\varepsilon)}{2-\alpha} \lim_{m\to\infty}\EE\left[\frac{1}{m}\log \frac{G_{k_m}}{G'_{k'_m}}\right] 
+ \frac{\varepsilon}{2-\alpha} \gamma(\alpha-\varepsilon). 
$$
Therefore, we obtain
\begin{eqnarray}
\label{on y arrive}
\frac{\gamma(\alpha) - \gamma(\alpha-\varepsilon)}{\varepsilon}
&=&\frac{\sigma(\alpha)-\sigma(\alpha-\varepsilon)}{\varepsilon}\dfrac{\gamma(\alpha)}{\sigma(\alpha)}\  +\  \frac{1}{2-\alpha} \gamma(\alpha-\varepsilon)\nonumber\\
&&\quad +\ \frac{\sigma(\alpha-\varepsilon)}{(2-\alpha)}\ \frac{1}{\varepsilon}\  \lim_{m\to\infty}\EE\left[\frac{1}{m}\log \frac{G_{k_m}}{G'_{k'_m}}\right] \, .
\end{eqnarray}
An easy computation shows that 
\begin{equation}
\label{1er_terme}
\lim_{\varepsilon\to 0} \left(\frac{\sigma(\alpha)-\sigma(\alpha-\varepsilon)}{\varepsilon}\dfrac{\gamma(\alpha)}{\sigma(\alpha)} + \frac{1}{2-\alpha} \gamma(\alpha-\varepsilon)\right)
=\gamma(\alpha)\frac{1+2\alpha-2\alpha^2}{(2\alpha-1)(\alpha^2-\alpha+1)}\, .
\end{equation}
Let us now turn to the last term on the right handside of \eqref{on y arrive}. 
For all $i\ge 1$, let us consider the $i$-th piece which is different in $Y$ and $Y'$. 
We denote by $m_i$ the number of pieces which have been seen before the $i$-th change. 
Let $a_i\egdef G_{k_{m_i}-1}$ and $b_i\egdef G_{k_{m_i}}$ be the labels along $Y$ of the vertices of the edge preceding the $i$-th change. 
Similarly, let $a'_i\egdef G_{k'_{m_i}-1}$ and $b'_i\egdef G'_{k'_{m_i}}$ be the corresponding labels in $Y'$. 
Using the fact that $i/m_i\to \varepsilon$ a.s., we get 
\begin{equation}
\label{2nd_terme}
\lim_{m\to\infty}\EE\left[\frac{1}{m}\log \frac{G_{k_m}}{G'_{k'_m}}\right]
= \varepsilon \lim_{i\to\infty} \EE\left[\frac{1}{i}\log \frac{a_i}{a'_i}\right]
\end{equation}
Between the $i$-th change and the edge $(a_{i+1}, b_{i+1})$, $Y$ and $Y'$ share a common part 
$Y^{i\to i+1}$ (see Figure~\ref{joining}). 
We set $n_i\egdef n(Y^{i\to i+1})$ and $d_i\egdef d(Y^{i\to i+1})$, where $n(\cdot)$ and $d(\cdot)$ are the notations introduced in the proof of Lemma~\ref{compare}. 
We then have the following induction.
$$
\frac{a_{i+1}}{a'_{i+1}} 
= \frac{d_ib_i+n_i(a_i+b_i)}{d_i(a'_i+b'_i)+n_ia'_i}
= \frac{a_{i}}{a'_{i}} \ 
\frac{\dfrac{b_i}{a_i}+\dfrac{n_i}{d_i}\left(1+\dfrac{b_i}{a_i}\right)}{1+\dfrac{b'_i}{a'_i}+\dfrac{n_i}{d_i}}\, .
$$
This yields
$$
\EE\left[\frac{1}{i}\log \frac{a_i}{a'_i}\right]
= \frac{1}{i} \sum_{j=1}^{i-1} \left(
\EE\left[ \log \left(\dfrac{b_j}{a_j}+\dfrac{n_j}{d_j}\left(1+\dfrac{b_j}{a_j}\right)\right)\right]
- \EE\left[ \log \left( 1+\dfrac{b'_j}{a'_j}+\dfrac{n_j}{d_j}\right)\right] \right)\, .
$$

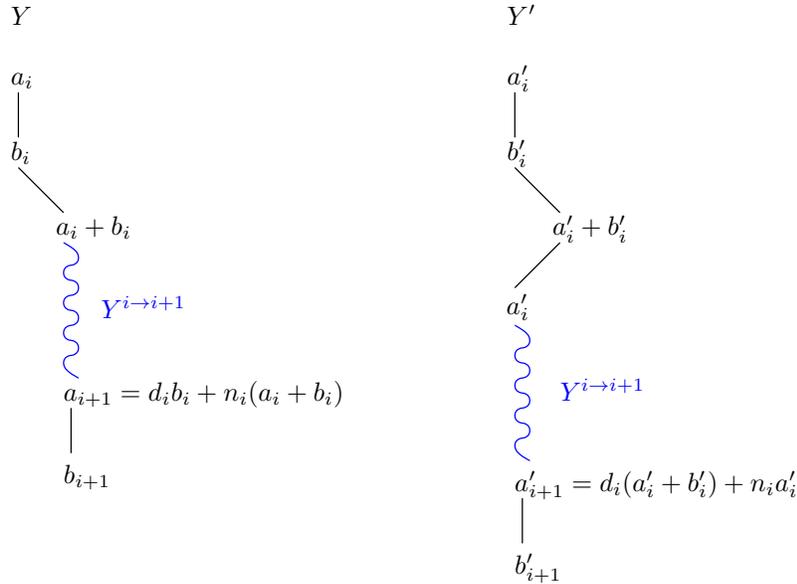
\begin{figure}[h]
	\begin{center}
	\input{joining.pstex_t}
	\end{center}
\caption{The paths $Y$ and $Y'$ between the $i$-th and $(i+1)$-th change.}
\label{joining}
\end{figure} 

Observe that $b_j/a_j$ and $n_j/d_j$ are independent. 
For all $j$, $n_j/d_j$ has a probability distribution $\nu_\alpha^\varepsilon$ which only depends on $\alpha$ and $\varepsilon$. 
Moreover, we know that $b_j/a_j$ converges in law to $\nu_\alpha$ (see Remark~\ref{remark nu_alpha}).
It follows that 
$$
\EE\left[ \log \left(\dfrac{b_j}{a_j}+\dfrac{n_j}{d_j}\left(1+\dfrac{b_j}{a_j}\right)\right)\right]
\tend{j}{\infty} \int \int 
\log \left(x+y\left(1+x\right)\right) d\nu_\alpha(x)d\nu_\alpha^\varepsilon(y)\, .
$$
Similarly, 
$$
\EE\left[ \log \left( 1+\dfrac{b'_j}{a'_j}+\dfrac{n_j}{d_j}\right)\right]
\tend{j}{\infty} \int \int 
\log \left(1+x+y\right) d\nu_\alpha(x)d\nu_\alpha^\varepsilon(y)\, .
$$
We thus obtain
$$
\EE\left[\frac{1}{i}\log \frac{a_i}{a'_i}\right] \tend{i}{\infty}\int \int 
\log \frac{x+y(1+x)}{1+x+y} d\nu_\alpha(x)d\nu_\alpha^\varepsilon(y)\, .
$$
The probability $\nu_\alpha^\varepsilon$ is the distribution of $n(Y^{1\to2})/d(Y^{1\to2})$. When $\varepsilon\to 0$, the length of the common part $Y^{1\to2}$ goes to infinity almost surely.
\begin{lemma}
\label{ouf}
Let $Y$ follow the probability distribution $\PP_\alpha$. Then 
$n(Y_1\dots Y_k)/d(Y_1\dots Y_k)$ has almost surely a limit as $k\to\infty$, which follows the probability distribution $\nu_\alpha$.
\end{lemma}
This lemma ensures that the preceding integral goes, as $\varepsilon\to 0$, to 
$$
\int \int \log \frac{x+y(1+x)}{1+x+y} d\nu_\alpha(x)d\nu_\alpha(y)\, .
$$
Together with \eqref{on y arrive}, \eqref{1er_terme} and \eqref{2nd_terme}, this achieves the proof. 
\end{proof}

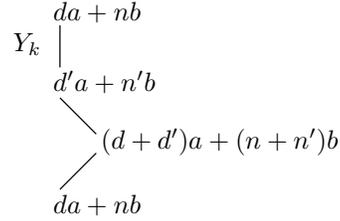
\begin{figure}[h]
	\begin{center}
	\input{nsurd.pstex_t}
	\end{center}
\caption{Labels on the path $Y$.}
\label{n_et_d}
\end{figure} 

\begin{proof}[Proof of Lemma~\ref{ouf}]
We decompose the path $Y$ into pieces $R$ or $RL$. Recall that, under $\PP_\alpha$, each piece appears independently with probability $\alpha$ for $R$ and $1-\alpha$ for $RL$. 
To each piece, we associate a real interval: Suppose the piece ends in $Y_k$. 
Then the bounds of the corresponding interval are defined as 
$n(Y_3\dots Y_{k-1})/d(Y_3\dots Y_{k-1})$ and $n(Y_3\dots Y_{k})/d(Y_3\dots Y_{k})$. 
Observe that the first interval is $[0, 1]$ with probability $1-\alpha$ and $[1, \infty]$ with probability $\alpha$. 
If $n/d=n(Y_3\dots Y_{k-1})/d(Y_3\dots Y_{k-1})$ and $n'/d'=n(Y_3\dots Y_{k})/d(Y_3\dots Y_{k})$ are the bounds of the interval associated to the $j$-th piece, then the bounds of the interval associated to the $(j+1)$-th piece are either $n'/d'$ and $(n+n')/(d+d')$ with probability $\alpha$, or $(n+n')/(d+d')$ and $n/d$ with probability $1-\alpha$ (see Figure~\ref{n_et_d}). 
Therefore, we get a decreasing sequence of Stern-Brocot intervals converging to a point following the probability distribution $\nu_\alpha$.
\end{proof}

\section{Link between $\nu_{\alpha}$ and Furstenberg's invariant measure}
In the linear case, Furstenberg's formula gives 
$$
\frac{1}{n} \log |F_n| \tend{n}{\infty} \gamma_p = \int \left( p \log\frac{\|xA\|}{\|x\|} + (1-p)\log\frac{\|xB\|}{\|x\|} \right)d\nu_f(\overline{x})\, ,
$$
where $A$ and $B$ are the matrices given in \eqref{matrices} and 
$\nu_f$ is the invariant measure on the set $P(\RR^2)$ of directions in the plane for the random walk that sends $\overline{x}$ to $\overline{x}A$ with probability $p$ and to $\overline{x}B$ with probability $1-p$. (In the above formula, $x$ stands for any nonzero vector with direction $\overline{x}$.)

Directions $\overline{x}$ can be parametrized using slopes $\overline{x} = (1, m)$ with $m\in (-\infty, \infty]$. 
Therefore, in this context, $\nu_f$ is the probability distribution on $(-\infty, \infty]$ such that, for any non-negative measurable function $g$,
$$
\int g \, d\nu_f = 
\int \left\{ p g\left(1+\frac{1}{m}\right) + (1-p) g\left(1-\frac{1}{m}\right)\right\} \, d\nu_f(m) \, .
$$
Observe that, in view of the particular form of the matrices $A$ and $B$ and the fact that $\nu_f$ is invariant, Furstenberg's formula reduces to 
\begin{eqnarray*}
\gamma_p 
& = &\int_{-\infty}^{+\infty} \log |m| \, d\nu_f(m) + 
\int \left( p \log\frac{\|(1, 1+\frac{1}{m})\|}{\|(1,m)\|} + (1-p)\log\frac{\|(1, 1-\frac{1}{m})\|}{\|(1,m)\|} \right)d\nu_f(m)\nonumber\\
& = &\int_{-\infty}^{+\infty} \log |m| \, d\nu_f(m)\, .
\end{eqnarray*}
(It is worth remarking that this simplification is always valid when dealing with linear recurrence equations.)

The difficult part is to identify the invariant measure $\nu_f$. 
This was done by Viswanath in the case $p=1/2$, and we could note that 
$\nu_f(\, .\, |\RR^+) = \nu_\alpha$ for $\alpha = \alpha(1/2) = \phi^{-1}$. 
This observation, together with the fact that the above equation looks very similar to our formula for $\gamma_p$, made us suspect a relationship between $\nu_\alpha$ and $\nu_f$. 

\begin{figure}[h]
	\begin{center}
	\input{mesure_nuf.pstex_t}
	\end{center}
\caption{The measure $\nu_f$ on Stern-Brocot intervals of rank 0, 1, 2.
First assign mass $m^- \egdef \dfrac{(1-\alpha)^2}{\alpha^2-\alpha+1}$ to $]-\infty,0]$ and $m^+\egdef\dfrac{\alpha}{\alpha^2-\alpha+1}$ to $[0,\infty]$. 
Once $\nu_\alpha$ is defined on some Stern-Brocot interval of rank $r$ in $\RR^+$, a proportion $\alpha$ of its mass is given to the left (respectively right) subinterval of rank $r+1$ when $r$ is odd (respectively even). In $\RR^-$, exchange $\alpha$ and $1-\alpha$.
}
\label{nuf}
\end{figure}
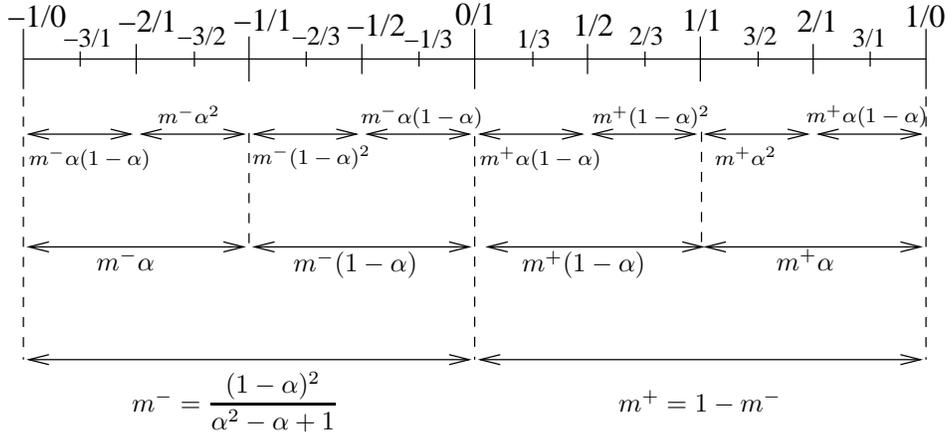

We came to the following heuristics: In the tree $\mathbf{T}$, all edges whose labels have opposite signs correspond to a step which will appear at the end of a deleted pattern $RLL$. 
Therefore, we expect $\nu_f(\RR^-)$ to be equal to the frequency of deletions 
$\dfrac{1-\sigma}{3} = \dfrac{(1-\alpha)^2}{\alpha^2-\alpha+1}$. 
Moreover, all edges whose labels have same signs can be seen as belonging to a tree $\mathbf{R}$. 
Therefore, we expect $\nu_f(\, .\, |\RR^+)$ to be equal to $\nu_\alpha$ (as in the case $p=1/2$).
The invariance property of $\nu_f$ for the indicator function of a Stern-Brocot interval $[a,b]$ with $1\le a<b$ now yields
$$ 
 \nu_f ([a,b]) = p\, \nu_f\left(\left[\dfrac{1}{b-1},\dfrac{1}{a-1}\right]\right)+ (1-p)\,\nu_f\left(\left[\dfrac{-1}{a-1},\dfrac{-1}{b-1}\right]\right).
$$
Observe that $\left[\frac{1}{b-1},\frac{1}{a-1}\right]$ and $\left[\frac{-1}{a-1},\frac{-1}{b-1}\right]$ are also Stern-Brocot intervals, the latter lying in $\RR^-$.
This equation is thus enough to get the measure of all Stern-Brocot intervals in $\RR^-$. We obtain the measure described in Figure~\ref{nuf}. 

Conversely, we easily check that this measure satisfies the invariance property.

\section{Further developments and open questions}
\subsection{Extension to Viswanath's setting}
The sequence $(t_n)$ studied by Viswanath in~\cite{viswanath2000} is defined by $t_1\egdef t_2\egdef1$, and 
$$ 
t_n = \pm t_{n-1}\pm t_{n-2}\quad (m\ge2),
$$
where each $\pm$ sign is independent and takes value $+$ or $-$ with probability $1/2$. A natural extension of this setting would be to choose the signs independently both with probability $p$ for $+$ and $1-p$ for $-$.
This definition is equivalent to ours only in the case $p=1/2$, which explains why the graph of the Lyapunov exponent drawn on Figure~5 in~\cite{viswanath2000} is different from our Figure~\ref{graphs_of_gamma}. In Viswanath's setting, no explicit formula is known to compute the Lyapunov exponent. Can our method be extended to this setting ?

\subsection{Random Fibonacci sequences with multiplicative coefficient} 
Consider the generalization of the random Fibonacci sequence $(F_{n})_{n}$ defined by
$F_{n+1}=\lambda F_{n}\pm F_{n-1}$ for a fixed parameter $\lambda$. 

For special values of this parameter, namely $\lambda_{k} = 2\cos(\pi/k)$, $k\ge 3$, we hope to generalize our method. 
(Observe that the present study corresponds to $k=3$). 
We expect the reduced sequences to be obtained by removing patterns $RL^{k-1}$, which should give rise to the study of a Markov chain of order $k-2$. 
The correspondance between random Fibonacci sequences and continued
fractions extends to these new sequences by considering {\em Rosen
continued fractions}, introduced by Rosen in \cite{Rosen1954}. A Rosen continued
fraction expansion of a real number $x$ is a continued fraction in which
partial quotients belong to $\lambda_k \ZZ^*$ (instead of $\ZZ_+^*$ as in the
classical case). The values $\lambda_k$ are known to be the only ones for
which the corresponding M\"obius
group generated by the transformations $z\longmapsto z-\lambda_k$ and
$z\longmapsto 1/\overline{z}$ acts discontinuously on $\HH$ \cite{GHJ1989,Magnus1975}. 

\subsection{Average growth rate} 
We announced in Remark~\ref{moyenne} that we were able to give an explicit expression for the limit of ${\dfrac{1}{n}\log(\EE({\widetilde F}_n))}$. 
Strangely enough, the similar question in the linear case seems to be more difficult. 
We do not know yet whether the combinatorial methods used in~\cite{rittaud2006} can be extended to this case.

\subsection{Variation properties}
We could expect the formula \eqref{derivee_gamma} to give the value of the derivative when $p\to0$ in the linear case (or when $p\to1/3$ in the non-linear case). Unfortunately an 
indeterminate form appears when $\alpha\to 1/2$. Is there another way to compute this value?

On figure~\ref{graphs_of_gamma}, obtained by numerical estimations of the integral, the functions $p\mapsto \gamma_p$ and $p\mapsto \widetilde\gamma_p$ seem to be convex. In~\cite{volkmer1999}, Volkmer gives a sufficient condition for the convexity of the Lyapunov exponent to hold. But this condition is easy to check only in the case of nonnegative matrices, and we do not know whether it applies in our cases.


\bibliography{rf.bib}

\end{document}

%% file: mesure.pstex_t
\begin{picture}(0,0)%
\includegraphics{mesure.pstex}%
\end{picture}%
\setlength{\unitlength}{4144sp}%
\begingroup\makeatletter\ifx\SetFigFont\undefined%
\gdef\SetFigFont#1#2#3#4#5{%
  \reset@font\fontsize{#1}{#2pt}%
  \fontfamily{#3}\fontseries{#4}\fontshape{#5}%
  \selectfont}%
\fi\endgroup%
\begin{picture}(5706,2419)(-149,-6830)
\put(766,-5416){\makebox(0,0)[lb]{\smash{{\SetFigFont{8}{9.6}{\rmdefault}{\mddefault}{\updefault}{\color[rgb]{0,0,0}$\alpha^2(1-\alpha)$}%
}}}}
\put(2746,-5416){\makebox(0,0)[lb]{\smash{{\SetFigFont{8}{9.6}{\rmdefault}{\mddefault}{\updefault}{\color[rgb]{0,0,0}$\alpha^2(1-\alpha)$}%
}}}}
\put(4816,-5416){\makebox(0,0)[lb]{\smash{{\SetFigFont{8}{9.6}{\rmdefault}{\mddefault}{\updefault}{\color[rgb]{0,0,0}$\alpha^2(1-\alpha)$}%
}}}}
\put(1441,-5416){\makebox(0,0)[lb]{\smash{{\SetFigFont{8}{9.6}{\rmdefault}{\mddefault}{\updefault}{\color[rgb]{0,0,0}$(1-\alpha)^3$}%
}}}}
\put(3646,-5416){\makebox(0,0)[lb]{\smash{{\SetFigFont{8}{9.6}{\rmdefault}{\mddefault}{\updefault}{\color[rgb]{0,0,0}$\alpha^3$}%
}}}}
\put(2116,-5416){\makebox(0,0)[lb]{\smash{{\SetFigFont{8}{9.6}{\rmdefault}{\mddefault}{\updefault}{\color[rgb]{0,0,0}$\alpha(1-\alpha)^2$}%
}}}}
\put(1666,-6091){\makebox(0,0)[lb]{\smash{{\SetFigFont{10}{12.0}{\rmdefault}{\mddefault}{\updefault}{\color[rgb]{0,0,0}$(1-\alpha)^2$}%
}}}}
\put(3286,-6091){\makebox(0,0)[lb]{\smash{{\SetFigFont{10}{12.0}{\rmdefault}{\mddefault}{\updefault}{\color[rgb]{0,0,0}$\alpha^2$}%
}}}}
\put(4366,-6091){\makebox(0,0)[lb]{\smash{{\SetFigFont{10}{12.0}{\rmdefault}{\mddefault}{\updefault}{\color[rgb]{0,0,0}$\alpha(1-\alpha)$}%
}}}}
\put( 46,-5416){\makebox(0,0)[lb]{\smash{{\SetFigFont{8}{9.6}{\rmdefault}{\mddefault}{\updefault}{\color[rgb]{0,0,0}$\alpha(1-\alpha)^2$}%
}}}}
\put(4141,-5416){\makebox(0,0)[lb]{\smash{{\SetFigFont{8}{9.6}{\rmdefault}{\mddefault}{\updefault}{\color[rgb]{0,0,0}$\alpha(1-\alpha)^2$}%
}}}}
\put(3916,-6766){\makebox(0,0)[lb]{\smash{{\SetFigFont{10}{12.0}{\rmdefault}{\mddefault}{\updefault}{\color[rgb]{0,0,0}$\alpha$}%
}}}}
\put(1036,-6766){\makebox(0,0)[lb]{\smash{{\SetFigFont{10}{12.0}{\rmdefault}{\mddefault}{\updefault}{\color[rgb]{0,0,0}$1-\alpha$}%
}}}}
\put(316,-6091){\makebox(0,0)[lb]{\smash{{\SetFigFont{10}{12.0}{\rmdefault}{\mddefault}{\updefault}{\color[rgb]{0,0,0}$\alpha(1-\alpha)$}%
}}}}
\end{picture}%

%% file: gammas.pstex_t
\begin{picture}(0,0)%
\includegraphics{gammas.pstex}%
\end{picture}%
\setlength{\unitlength}{4144sp}%
\begingroup\makeatletter\ifx\SetFigFont\undefined%
\gdef\SetFigFont#1#2#3#4#5{%
  \reset@font\fontsize{#1}{#2pt}%
  \fontfamily{#3}\fontseries{#4}\fontshape{#5}%
  \selectfont}%
\fi\endgroup%
\begin{picture}(6017,1789)(-356,608)
\put(865,2250){\makebox(0,0)[lb]{\smash{{\SetFigFont{10}{12.0}{\rmdefault}{\mddefault}{\updefault}{\color[rgb]{0,0,0}$\gamma_p$}%
}}}}
\put(865,2004){\makebox(0,0)[lb]{\smash{{\SetFigFont{10}{12.0}{\rmdefault}{\mddefault}{\updefault}{\color[rgb]{0,0,0}$\widetilde\gamma_p$}%
}}}}
\put(3316,669){\makebox(0,0)[lb]{\smash{{\SetFigFont{10}{12.0}{\rmdefault}{\mddefault}{\updefault}{\color[rgb]{0,0,0}0}%
}}}}
\put(5561,672){\makebox(0,0)[lb]{\smash{{\SetFigFont{10}{12.0}{\rmdefault}{\mddefault}{\updefault}{\color[rgb]{0,0,0}$\alpha$}%
}}}}
\put(3033,1766){\makebox(0,0)[lb]{\smash{{\SetFigFont{10}{12.0}{\rmdefault}{\mddefault}{\updefault}{\color[rgb]{0,0,0}$\ln \phi$}%
}}}}
\put(5313,672){\makebox(0,0)[lb]{\smash{{\SetFigFont{10}{12.0}{\rmdefault}{\mddefault}{\updefault}{\color[rgb]{0,0,0}1}%
}}}}
\put(4251,686){\makebox(0,0)[lb]{\smash{{\SetFigFont{10}{12.0}{\rmdefault}{\mddefault}{\updefault}{\color[rgb]{0,0,0}$1/2$}%
}}}}
\put(-58,669){\makebox(0,0)[lb]{\smash{{\SetFigFont{10}{12.0}{\rmdefault}{\mddefault}{\updefault}{\color[rgb]{0,0,0}0}%
}}}}
\put(2187,672){\makebox(0,0)[lb]{\smash{{\SetFigFont{10}{12.0}{\rmdefault}{\mddefault}{\updefault}{\color[rgb]{0,0,0}$p$}%
}}}}
\put(-341,1766){\makebox(0,0)[lb]{\smash{{\SetFigFont{10}{12.0}{\rmdefault}{\mddefault}{\updefault}{\color[rgb]{0,0,0}$\ln \phi$}%
}}}}
\put(1942,657){\makebox(0,0)[lb]{\smash{{\SetFigFont{10}{12.0}{\rmdefault}{\mddefault}{\updefault}{\color[rgb]{0,0,0}1}%
}}}}
\put(3016,2038){\makebox(0,0)[lb]{\smash{{\SetFigFont{10}{12.0}{\rmdefault}{\mddefault}{\updefault}{\color[rgb]{0,0,0}$\gamma(\alpha)$}%
}}}}
\put(565,657){\makebox(0,0)[lb]{\smash{{\SetFigFont{10}{12.0}{\rmdefault}{\mddefault}{\updefault}{\color[rgb]{0,0,0}1/3}%
}}}}
\put(932,657){\makebox(0,0)[lb]{\smash{{\SetFigFont{10}{12.0}{\rmdefault}{\mddefault}{\updefault}{\color[rgb]{0,0,0}1/2}%
}}}}
\end{picture}%

%% file: reduction.pstex_t
\begin{picture}(0,0)%
\includegraphics{reduction.pstex}%
\end{picture}%
\setlength{\unitlength}{2901sp}%
\begingroup\makeatletter\ifx\SetFigFont\undefined%
\gdef\SetFigFont#1#2#3#4#5{%
  \reset@font\fontsize{#1}{#2pt}%
  \fontfamily{#3}\fontseries{#4}\fontshape{#5}%
  \selectfont}%
\fi\endgroup%
\begin{picture}(4316,3242)(1837,-5673)
\put(1852,-5591){\makebox(0,0)[lb]{\smash{{\SetFigFont{8}{9.6}{\rmdefault}{\mddefault}{\itdefault}{\color[rgb]{0,0,0}$-(b+a)$}%
}}}}
\put(2868,-5600){\makebox(0,0)[lb]{\smash{{\SetFigFont{8}{9.6}{\rmdefault}{\mddefault}{\itdefault}{\color[rgb]{0,0,0}$-(b-a)$}%
}}}}
\put(3061,-3211){\makebox(0,0)[lb]{\smash{{\SetFigFont{8}{9.6}{\rmdefault}{\mddefault}{\itdefault}{\color[rgb]{0,0,0}$b$}%
}}}}
\put(2568,-3774){\makebox(0,0)[lb]{\smash{{\SetFigFont{8}{9.6}{\rmdefault}{\mddefault}{\itdefault}{\color[rgb]{0,0,0}$b-a$}%
}}}}
\put(3579,-3771){\makebox(0,0)[lb]{\smash{{\SetFigFont{8}{9.6}{\rmdefault}{\mddefault}{\itdefault}{\color[rgb]{0,0,0}$b+a$}%
}}}}
\put(3065,-2578){\makebox(0,0)[lb]{\smash{{\SetFigFont{8}{9.6}{\rmdefault}{\mddefault}{\itdefault}{\color[rgb]{0,0,0}$a$}%
}}}}
\put(3065,-4378){\makebox(0,0)[lb]{\smash{{\SetFigFont{8}{9.6}{\rmdefault}{\mddefault}{\itdefault}{\color[rgb]{0,0,0}$a$}%
}}}}
\put(2487,-4977){\makebox(0,0)[lb]{\smash{{\SetFigFont{8}{9.6}{\rmdefault}{\mddefault}{\itdefault}{\color[rgb]{0,0,0}$-b$}%
}}}}
\put(3442,-3418){\makebox(0,0)[lb]{\smash{{\SetFigFont{8}{9.6}{\rmdefault}{\mddefault}{\updefault}{\color[rgb]{0,0,1}$R$}%
}}}}
\put(3482,-4111){\makebox(0,0)[lb]{\smash{{\SetFigFont{8}{9.6}{\rmdefault}{\mddefault}{\updefault}{\color[rgb]{0,0,1}$L$}%
}}}}
\put(2922,-4705){\makebox(0,0)[lb]{\smash{{\SetFigFont{8}{9.6}{\rmdefault}{\mddefault}{\updefault}{\color[rgb]{0,0,1}$L$}%
}}}}
\put(5616,-3211){\makebox(0,0)[lb]{\smash{{\SetFigFont{8}{9.6}{\rmdefault}{\mddefault}{\itdefault}{\color[rgb]{0,0,0}$b$}%
}}}}
\put(6134,-3771){\makebox(0,0)[lb]{\smash{{\SetFigFont{8}{9.6}{\rmdefault}{\mddefault}{\itdefault}{\color[rgb]{0,0,0}$b+a$}%
}}}}
\put(5620,-2578){\makebox(0,0)[lb]{\smash{{\SetFigFont{8}{9.6}{\rmdefault}{\mddefault}{\itdefault}{\color[rgb]{0,0,0}$a$}%
}}}}
\put(5620,-4378){\makebox(0,0)[lb]{\smash{{\SetFigFont{8}{9.6}{\rmdefault}{\mddefault}{\itdefault}{\color[rgb]{0,0,0}$a$}%
}}}}
\put(5042,-4977){\makebox(0,0)[lb]{\smash{{\SetFigFont{8}{9.6}{\rmdefault}{\mddefault}{\itdefault}{\color[rgb]{0,0,0}$b$}%
}}}}
\put(6013,-3418){\makebox(0,0)[lb]{\smash{{\SetFigFont{8}{9.6}{\rmdefault}{\mddefault}{\updefault}{\color[rgb]{0,0,1}$R$}%
}}}}
\put(6046,-4111){\makebox(0,0)[lb]{\smash{{\SetFigFont{8}{9.6}{\rmdefault}{\mddefault}{\updefault}{\color[rgb]{0,0,1}$L$}%
}}}}
\put(5473,-4705){\makebox(0,0)[lb]{\smash{{\SetFigFont{8}{9.6}{\rmdefault}{\mddefault}{\updefault}{\color[rgb]{0,0,1}$L$}%
}}}}
\end{picture}%

%% file: probaR.pstex_t
\begin{picture}(0,0)%
\includegraphics{probaR.pstex}%
\end{picture}%
\setlength{\unitlength}{4144sp}%
\begingroup\makeatletter\ifx\SetFigFont\undefined%
\gdef\SetFigFont#1#2#3#4#5{%
  \reset@font\fontsize{#1}{#2pt}%
  \fontfamily{#3}\fontseries{#4}\fontshape{#5}%
  \selectfont}%
\fi\endgroup%
\begin{picture}(2728,1360)(2368,-3838)
\put(3061,-3211){\makebox(0,0)[lb]{\smash{{\SetFigFont{12}{14.4}{\rmdefault}{\mddefault}{\itdefault}{\color[rgb]{0,0,0}$b$}%
}}}}
\put(2568,-3774){\makebox(0,0)[lb]{\smash{{\SetFigFont{12}{14.4}{\rmdefault}{\mddefault}{\itdefault}{\color[rgb]{0,0,0}$b-a$}%
}}}}
\put(3433,-3421){\makebox(0,0)[lb]{\smash{{\SetFigFont{12}{14.4}{\rmdefault}{\mddefault}{\updefault}{\color[rgb]{0,0,0}$\alpha$}%
}}}}
\put(2383,-3421){\makebox(0,0)[lb]{\smash{{\SetFigFont{12}{14.4}{\rmdefault}{\mddefault}{\updefault}{\color[rgb]{0,0,0}$1-\alpha$}%
}}}}
\put(3579,-3771){\makebox(0,0)[lb]{\smash{{\SetFigFont{12}{14.4}{\rmdefault}{\mddefault}{\itdefault}{\color[rgb]{0,0,0}$b+a$}%
}}}}
\put(2626,-2625){\makebox(0,0)[lb]{\smash{{\SetFigFont{12}{14.4}{\rmdefault}{\mddefault}{\itdefault}{\color[rgb]{0,0,0}$a$}%
}}}}
\put(4990,-2625){\makebox(0,0)[lb]{\smash{{\SetFigFont{12}{14.4}{\rmdefault}{\mddefault}{\itdefault}{\color[rgb]{0,0,0}$a$}%
}}}}
\put(4600,-3211){\makebox(0,0)[lb]{\smash{{\SetFigFont{12}{14.4}{\rmdefault}{\mddefault}{\itdefault}{\color[rgb]{0,0,0}$b$}%
}}}}
\put(5011,-3771){\makebox(0,0)[lb]{\smash{{\SetFigFont{12}{14.4}{\rmdefault}{\mddefault}{\itdefault}{\color[rgb]{0,0,0}$b+a$}%
}}}}
\end{picture}%

%% file: ex_fraction.pstex_t
\begin{picture}(0,0)%
\includegraphics{ex_fraction.pstex}%
\end{picture}%
\setlength{\unitlength}{4144sp}%
\begingroup\makeatletter\ifx\SetFigFont\undefined%
\gdef\SetFigFont#1#2#3#4#5{%
  \reset@font\fontsize{#1}{#2pt}%
  \fontfamily{#3}\fontseries{#4}\fontshape{#5}%
  \selectfont}%
\fi\endgroup%
\begin{picture}(2368,4736)(6916,-9102)
\put(9219,-8997){\makebox(0,0)[lb]{\smash{{\SetFigFont{12}{14.4}{\rmdefault}{\mddefault}{\updefault}{\color[rgb]{0,0,0}$a_1 = 1$}%
}}}}
\put(9219,-8324){\makebox(0,0)[lb]{\smash{{\SetFigFont{12}{14.4}{\rmdefault}{\mddefault}{\updefault}{\color[rgb]{0,0,0}$a_2 = 2$}%
}}}}
\put(9219,-7591){\makebox(0,0)[lb]{\smash{{\SetFigFont{12}{14.4}{\rmdefault}{\mddefault}{\updefault}{\color[rgb]{0,0,0}$a_3 = 1$}%
}}}}
\put(9219,-5917){\makebox(0,0)[lb]{\smash{{\SetFigFont{12}{14.4}{\rmdefault}{\mddefault}{\updefault}{\color[rgb]{0,0,0}$a_4 = 6$}%
}}}}
\end{picture}%

%% file: evolvQ.pstex_t
\begin{picture}(0,0)%
\includegraphics{evolvQ.pstex}%
\end{picture}%
\setlength{\unitlength}{3315sp}%
\begingroup\makeatletter\ifx\SetFigFont\undefined%
\gdef\SetFigFont#1#2#3#4#5{%
  \reset@font\fontsize{#1}{#2pt}%
  \fontfamily{#3}\fontseries{#4}\fontshape{#5}%
  \selectfont}%
\fi\endgroup%
\begin{picture}(7827,3495)(-101,-4013)
\put(4717,-1810){\makebox(0,0)[lb]{\smash{{\SetFigFont{10}{12.0}{\rmdefault}{\mddefault}{\updefault}{\color[rgb]{0,0,0}$x$}%
}}}}
\put(7144,-2158){\makebox(0,0)[lb]{\smash{{\SetFigFont{10}{12.0}{\rmdefault}{\mddefault}{\updefault}{\color[rgb]{0,0,0}$q_{k-1}$}%
}}}}
\put(7474,-2570){\makebox(0,0)[lb]{\smash{{\SetFigFont{10}{12.0}{\rmdefault}{\mddefault}{\updefault}{\color[rgb]{0,0,0}$y$}%
}}}}
\put(7158,-2917){\makebox(0,0)[lb]{\smash{{\SetFigFont{10}{12.0}{\rmdefault}{\mddefault}{\updefault}{\color[rgb]{0,0,0}$q_{k}$}%
}}}}
\put(6572,-3293){\makebox(0,0)[lb]{\smash{{\SetFigFont{10}{12.0}{\rmdefault}{\mddefault}{\updefault}{\color[rgb]{0,0,0}$y-x$}%
}}}}
\put(1797,-1744){\makebox(0,0)[lb]{\smash{{\SetFigFont{10}{12.0}{\rmdefault}{\mddefault}{\updefault}{\color[rgb]{0,0,0}$x$}%
}}}}
\put(2252,-2092){\makebox(0,0)[lb]{\smash{{\SetFigFont{10}{12.0}{\rmdefault}{\mddefault}{\updefault}{\color[rgb]{0,0,0}$q_{k-1}$}%
}}}}
\put(3020,-2870){\makebox(0,0)[lb]{\smash{{\SetFigFont{10}{12.0}{\rmdefault}{\mddefault}{\updefault}{\color[rgb]{0,0,0}$q_{k}$}%
}}}}
\put(3154,-3345){\makebox(0,0)[lb]{\smash{{\SetFigFont{10}{12.0}{\rmdefault}{\mddefault}{\updefault}{\color[rgb]{0,0,0}$x+y$}%
}}}}
\put(2594,-2442){\makebox(0,0)[lb]{\smash{{\SetFigFont{10}{12.0}{\rmdefault}{\mddefault}{\updefault}{\color[rgb]{0,0,0}$y$}%
}}}}
\put(6446,-1788){\makebox(0,0)[lb]{\smash{{\SetFigFont{10}{12.0}{\rmdefault}{\mddefault}{\updefault}{\color[rgb]{0,0,0}$x$}%
}}}}
\put(787,-1908){\makebox(0,0)[lb]{\smash{{\SetFigFont{10}{12.0}{\rmdefault}{\mddefault}{\updefault}{\color[rgb]{0,0,0}$x$}%
}}}}
\put(609,-3362){\makebox(0,0)[lb]{\smash{{\SetFigFont{10}{12.0}{\rmdefault}{\mddefault}{\updefault}{\color[rgb]{0,0,0}$x+y$}%
}}}}
\put(463,-2917){\makebox(0,0)[lb]{\smash{{\SetFigFont{10}{12.0}{\rmdefault}{\mddefault}{\updefault}{\color[rgb]{0,0,0}$q_{k}$}%
}}}}
\put(-86,-2570){\makebox(0,0)[lb]{\smash{{\SetFigFont{10}{12.0}{\rmdefault}{\mddefault}{\updefault}{\color[rgb]{0,0,0}$y$}%
}}}}
\put(5172,-2158){\makebox(0,0)[lb]{\smash{{\SetFigFont{10}{12.0}{\rmdefault}{\mddefault}{\updefault}{\color[rgb]{0,0,0}$q_{k-1}$}%
}}}}
\put(5502,-2570){\makebox(0,0)[lb]{\smash{{\SetFigFont{10}{12.0}{\rmdefault}{\mddefault}{\updefault}{\color[rgb]{0,0,0}$y$}%
}}}}
\put(5186,-2917){\makebox(0,0)[lb]{\smash{{\SetFigFont{10}{12.0}{\rmdefault}{\mddefault}{\updefault}{\color[rgb]{0,0,0}$q_{k}$}%
}}}}
\put(4600,-3293){\makebox(0,0)[lb]{\smash{{\SetFigFont{10}{12.0}{\rmdefault}{\mddefault}{\updefault}{\color[rgb]{0,0,0}$y-x$}%
}}}}
\put(-44,-2158){\makebox(0,0)[lb]{\smash{{\SetFigFont{10}{12.0}{\rmdefault}{\mddefault}{\updefault}{\color[rgb]{0,0,0}$q_{k-1}$}%
}}}}
\end{picture}%

%% file: def_n_et_d.pstex_t
\begin{picture}(0,0)%
\includegraphics{def_n_et_d.pstex}%
\end{picture}%
\setlength{\unitlength}{4144sp}%
\begingroup\makeatletter\ifx\SetFigFont\undefined%
\gdef\SetFigFont#1#2#3#4#5{%
  \reset@font\fontsize{#1}{#2pt}%
  \fontfamily{#3}\fontseries{#4}\fontshape{#5}%
  \selectfont}%
\fi\endgroup%
\begin{picture}(1519,2118)(-198,-3136)
\put(676,-1951){\makebox(0,0)[lb]{\smash{{\SetFigFont{10}{12.0}{\rmdefault}{\mddefault}{\updefault}{\color[rgb]{0,0,0}$a+b$}%
}}}}
\put(451,-1546){\makebox(0,0)[lb]{\smash{{\SetFigFont{10}{12.0}{\rmdefault}{\mddefault}{\updefault}{\color[rgb]{0,0,0}$b$}%
}}}}
\put(451,-1141){\makebox(0,0)[lb]{\smash{{\SetFigFont{10}{12.0}{\rmdefault}{\mddefault}{\updefault}{\color[rgb]{0,0,0}$a$}%
}}}}
\put(856,-3076){\makebox(0,0)[lb]{\smash{{\SetFigFont{10}{12.0}{\rmdefault}{\mddefault}{\updefault}{\color[rgb]{0,0,0}$d(y)a+n(y)b$}%
}}}}
\put(1306,-1951){\makebox(0,0)[lb]{\smash{{\SetFigFont{10}{12.0}{\rmdefault}{\mddefault}{\updefault}{\color[rgb]{0,0,0}$d(R)=1$, $n(R)=1$}%
}}}}
\put(1306,-1546){\makebox(0,0)[lb]{\smash{{\SetFigFont{10}{12.0}{\rmdefault}{\mddefault}{\updefault}{\color[rgb]{0,0,0}$d(\emptyset)=0$, $n(\emptyset)=1$}%
}}}}
\put(-183,-2446){\makebox(0,0)[lb]{\smash{{\SetFigFont{10}{12.0}{\rmdefault}{\mddefault}{\updefault}{\color[rgb]{0,0,0}$y$}%
}}}}
\end{picture}%

%% file: joining.pstex_t
\begin{picture}(0,0)%
\includegraphics{joining.pstex}%
\end{picture}%
\setlength{\unitlength}{4144sp}%
\begingroup\makeatletter\ifx\SetFigFont\undefined%
\gdef\SetFigFont#1#2#3#4#5{%
  \reset@font\fontsize{#1}{#2pt}%
  \fontfamily{#3}\fontseries{#4}\fontshape{#5}%
  \selectfont}%
\fi\endgroup%
\begin{picture}(5032,3602)(707,-2633)
\put(811,389){\makebox(0,0)[lb]{\smash{{\SetFigFont{10}{12.0}{\rmdefault}{\mddefault}{\updefault}{\color[rgb]{0,0,0}$a_i$}%
}}}}
\put(811,-61){\makebox(0,0)[lb]{\smash{{\SetFigFont{10}{12.0}{\rmdefault}{\mddefault}{\updefault}{\color[rgb]{0,0,0}$b_i$}%
}}}}
\put(1081,-511){\makebox(0,0)[lb]{\smash{{\SetFigFont{10}{12.0}{\rmdefault}{\mddefault}{\updefault}{\color[rgb]{0,0,0}$a_i+b_i$}%
}}}}
\put(3781,389){\makebox(0,0)[lb]{\smash{{\SetFigFont{10}{12.0}{\rmdefault}{\mddefault}{\updefault}{\color[rgb]{0,0,0}$a'_i$}%
}}}}
\put(3781,-61){\makebox(0,0)[lb]{\smash{{\SetFigFont{10}{12.0}{\rmdefault}{\mddefault}{\updefault}{\color[rgb]{0,0,0}$b'_i$}%
}}}}
\put(4051,-511){\makebox(0,0)[lb]{\smash{{\SetFigFont{10}{12.0}{\rmdefault}{\mddefault}{\updefault}{\color[rgb]{0,0,0}$a'_i+b'_i$}%
}}}}
\put(811,749){\makebox(0,0)[lb]{\smash{{\SetFigFont{10}{12.0}{\rmdefault}{\mddefault}{\updefault}{\color[rgb]{0,0,0}$Y$}%
}}}}
\put(3781,749){\makebox(0,0)[lb]{\smash{{\SetFigFont{10}{12.0}{\rmdefault}{\mddefault}{\updefault}{\color[rgb]{0,0,0}$Y'$}%
}}}}
\put(3777,-980){\makebox(0,0)[lb]{\smash{{\SetFigFont{10}{12.0}{\rmdefault}{\mddefault}{\updefault}{\color[rgb]{0,0,0}$a'_i$}%
}}}}
\put(1126,-1996){\makebox(0,0)[lb]{\smash{{\SetFigFont{10}{12.0}{\rmdefault}{\mddefault}{\updefault}{\color[rgb]{0,0,0}$b_{i+1}$}%
}}}}
\put(3826,-2041){\makebox(0,0)[lb]{\smash{{\SetFigFont{10}{12.0}{\rmdefault}{\mddefault}{\updefault}{\color[rgb]{0,0,0}$a'_{i+1}=d_i(a'_i+b'_i)+n_ia'_i$}%
}}}}
\put(3826,-2536){\makebox(0,0)[lb]{\smash{{\SetFigFont{10}{12.0}{\rmdefault}{\mddefault}{\updefault}{\color[rgb]{0,0,0}$b'_{i+1}$}%
}}}}
\put(1126,-1501){\makebox(0,0)[lb]{\smash{{\SetFigFont{10}{12.0}{\rmdefault}{\mddefault}{\updefault}{\color[rgb]{0,0,0}$a_{i+1}=d_ib_i+n_i(a_i+b_i)$}%
}}}}
\put(1351,-1006){\makebox(0,0)[lb]{\smash{{\SetFigFont{10}{12.0}{\rmdefault}{\mddefault}{\updefault}{\color[rgb]{0,0,1}$Y^{i\to i+1}$}%
}}}}
\put(4096,-1501){\makebox(0,0)[lb]{\smash{{\SetFigFont{10}{12.0}{\rmdefault}{\mddefault}{\updefault}{\color[rgb]{0,0,1}$Y^{i\to i+1}$}%
}}}}
\end{picture}%

%% file: nsurd.pstex_t
\begin{picture}(0,0)%
\includegraphics{nsurd.pstex}%
\end{picture}%
\setlength{\unitlength}{4144sp}%
\begingroup\makeatletter\ifx\SetFigFont\undefined%
\gdef\SetFigFont#1#2#3#4#5{%
  \reset@font\fontsize{#1}{#2pt}%
  \fontfamily{#3}\fontseries{#4}\fontshape{#5}%
  \selectfont}%
\fi\endgroup%
\begin{picture}(556,1332)(273,-939)
\put(529,-151){\makebox(0,0)[lb]{\smash{{\SetFigFont{10}{12.0}{\rmdefault}{\mddefault}{\updefault}{\color[rgb]{0,0,0}$d'a+n'b$}%
}}}}
\put(814,-496){\makebox(0,0)[lb]{\smash{{\SetFigFont{10}{12.0}{\rmdefault}{\mddefault}{\updefault}{\color[rgb]{0,0,0}$(d+d')a+(n+n')b$}%
}}}}
\put(529,-884){\makebox(0,0)[lb]{\smash{{\SetFigFont{10}{12.0}{\rmdefault}{\mddefault}{\updefault}{\color[rgb]{0,0,0}$da+nb$}%
}}}}
\put(529,270){\makebox(0,0)[lb]{\smash{{\SetFigFont{10}{12.0}{\rmdefault}{\mddefault}{\updefault}{\color[rgb]{0,0,0}$da+nb$}%
}}}}
\put(288, 85){\makebox(0,0)[lb]{\smash{{\SetFigFont{10}{12.0}{\rmdefault}{\mddefault}{\updefault}{\color[rgb]{0,0,0}$Y_k$}%
}}}}
\end{picture}%

%% file: mesure_nuf.pstex_t
\begin{picture}(0,0)%
\includegraphics{mesure_nuf.pstex}%
\end{picture}%
\setlength{\unitlength}{4144sp}%
\begingroup\makeatletter\ifx\SetFigFont\undefined%
\gdef\SetFigFont#1#2#3#4#5{%
  \reset@font\fontsize{#1}{#2pt}%
  \fontfamily{#3}\fontseries{#4}\fontshape{#5}%
  \selectfont}%
\fi\endgroup%
\begin{picture}(5672,2526)(-109,-6957)
\put(2983,-6057){\makebox(0,0)[lb]{\smash{{\SetFigFont{10}{12.0}{\rmdefault}{\mddefault}{\updefault}{\color[rgb]{0,0,0}$m^+(1-\alpha)$}%
}}}}
\put(1617,-6057){\makebox(0,0)[lb]{\smash{{\SetFigFont{10}{12.0}{\rmdefault}{\mddefault}{\updefault}{\color[rgb]{0,0,0}$m^-(1-\alpha)$}%
}}}}
\put(438,-6051){\makebox(0,0)[lb]{\smash{{\SetFigFont{10}{12.0}{\rmdefault}{\mddefault}{\updefault}{\color[rgb]{0,0,0}$m^-\alpha$}%
}}}}
\put(2732,-5418){\makebox(0,0)[lb]{\smash{{\SetFigFont{8}{9.6}{\rmdefault}{\mddefault}{\updefault}{\color[rgb]{0,0,0}$m^+\alpha(1-\alpha)$}%
}}}}
\put(4506,-6051){\makebox(0,0)[lb]{\smash{{\SetFigFont{10}{12.0}{\rmdefault}{\mddefault}{\updefault}{\color[rgb]{0,0,0}$m^+\alpha$}%
}}}}
\put(807,-5171){\makebox(0,0)[lb]{\smash{{\SetFigFont{8}{9.6}{\rmdefault}{\mddefault}{\updefault}{\color[rgb]{0,0,0}$m^-\alpha^2$}%
}}}}
\put( 34,-5417){\makebox(0,0)[lb]{\smash{{\SetFigFont{8}{9.6}{\rmdefault}{\mddefault}{\updefault}{\color[rgb]{0,0,0}$m^-\alpha(1-\alpha)$}%
}}}}
\put(1373,-5411){\makebox(0,0)[lb]{\smash{{\SetFigFont{8}{9.6}{\rmdefault}{\mddefault}{\updefault}{\color[rgb]{0,0,0}$m^-(1-\alpha)^2$}%
}}}}
\put(3406,-5171){\makebox(0,0)[lb]{\smash{{\SetFigFont{8}{9.6}{\rmdefault}{\mddefault}{\updefault}{\color[rgb]{0,0,0}$m^+(1-\alpha)^2$}%
}}}}
\put(649,-6893){\makebox(0,0)[lb]{\smash{{\SetFigFont{10}{12.0}{\rmdefault}{\mddefault}{\updefault}{\color[rgb]{0,0,0}$m^-=\dfrac{(1-\alpha)^2}{\alpha^2-\alpha+1}$}%
}}}}
\put(3562,-6893){\makebox(0,0)[lb]{\smash{{\SetFigFont{10}{12.0}{\rmdefault}{\mddefault}{\updefault}{\color[rgb]{0,0,0}$m^+=1-m^-$}%
}}}}
\put(4139,-5417){\makebox(0,0)[lb]{\smash{{\SetFigFont{8}{9.6}{\rmdefault}{\mddefault}{\updefault}{\color[rgb]{0,0,0}$m^+\alpha^2$}%
}}}}
\put(4686,-5165){\makebox(0,0)[lb]{\smash{{\SetFigFont{8}{9.6}{\rmdefault}{\mddefault}{\updefault}{\color[rgb]{0,0,0}$m^+\alpha(1-\alpha)$}%
}}}}
\put(2022,-5169){\makebox(0,0)[lb]{\smash{{\SetFigFont{8}{9.6}{\rmdefault}{\mddefault}{\updefault}{\color[rgb]{0,0,0}$m^-\alpha(1-\alpha)$}%
}}}}
\end{picture}%

%% file: rf.bbl
\providecommand{\bysame}{\leavevmode\hbox to3em{\hrulefill}\thinspace}
\providecommand{\MR}{\relax\ifhmode\unskip\space\fi MR }
\providecommand{\MRhref}[2]{%
  \href{http://www.ams.org/mathscinet-getitem?mr=#1}{#2}
}
\providecommand{\href}[2]{#2}
\begin{thebibliography}{10}

\bibitem{bougerol1985}
Philippe Bougerol and Jean Lacroix, \emph{Products of random matrices with
  applications to {S}chr\"odinger operators}, Progress in Probability and
  Statistics, vol.~8, Birkh\"auser Boston Inc., Boston, MA, 1985.

\bibitem{chassaing1984}
Philippe Chassaing, G{\'e}rard Letac, and Marianne Mora, \emph{Brocot sequences
  and random walks in {${\rm SL}(2,{\bf R})$}}, Probability measures on groups,
  VII (Oberwolfach, 1983), Lecture Notes in Math., vol. 1064, Springer, Berlin,
  1984, pp.~36--48.

\bibitem{denjoy1938}
Arnaud Denjoy, \emph{Sur une fonction r{\'e}elle de {M}inkowski}, J. Math.
  Pures Appl. \textbf{17} (1938), 105--151.

\bibitem{furstenberg1963}
Harry Furstenberg, \emph{Noncommuting random products}, Trans. Amer. Math. Soc.
  \textbf{108} (1963), 377--428.

\bibitem{GHJ1989}
Frederick~M. Goodman, Pierre de~la Harpe, and Vaughan F.~R. Jones,
  \emph{Coxeter graphs and towers of algebras}, Mathematical Sciences Research
  Institute Publications, vol.~14, Springer-Verlag, New York, 1989.

\bibitem{Magnus1975}
Wilhelm Magnus, \emph{Two generator subgroups of {${\rm PSL}$} {$(2,\,C)$}},
  Nachr. Akad. Wiss. G\"ottingen Math.-Phys. Kl. II (1975), no.~7, 81--94.

\bibitem{peres1990}
Yuval Peres, \emph{Analytic dependence of {L}yapunov exponents on transition
  probabilities}, Lyapunov exponents (Oberwolfach, 1990), Lecture Notes in
  Math., vol. 1486, Springer, Berlin, 1991, pp.~64--80.

\bibitem{rittaud2006}
Beno\^it Rittaud, \emph{On the average growth of random {F}ibonacci sequences},
  to appear in Journal of Integer Sequences, 2006.

\bibitem{Rosen1954}
David Rosen, \emph{A class of continued fractions associated with certain
  properly discontinuous groups}, Duke Math. J. \textbf{21} (1954), 549--563.

\bibitem{viswanath2000}
Divakar Viswanath, \emph{Random {F}ibonacci sequences and the number
  {$1.13198824\ldots$}}, Math. Comp. \textbf{69} (2000), no.~231, 1131--1155.

\bibitem{volkmer1999}
Hans Volkmer, \emph{Convexity of the {L}yapunov exponent}, Linear Algebra Appl.
  \textbf{294} (1999), no.~1-3, 35--48.

\end{thebibliography}
